\documentclass[10pt]{article}
\usepackage{amsmath,amssymb}
\usepackage[dvips]{graphicx}
\usepackage{latexsym}
\usepackage{bm}
\def\reclame #1\par{\medbreak \noindent{\hspace*{5mm} \it#1}\par\medbreak}

\newcommand{\no}[1]{, No.\,#1}   
\newcommand{\vol}[1]{{\bf #1}}   

\newcommand{\refstyle}{\it} 

\def\endrefs{\end{enumerate}\label{LASTPAGE}}

\newcommand{\proof}{\par\addvspace{0.5\baselineskip}\par{\bf Proof.\ }}

\def\endproof{\quad $\Box$\par\addvspace{0.3\baselineskip}\par}

\newcommand{\refs}{\par\subsubsection*{{\small\bf REFERENCES}}
\par%
\begin{enumerate}}          
\numberwithin{equation}{section}
\renewcommand{\d}{\,{\rm d}}    
\newcommand{\e}{\mkern 1mu{\rm e}}               
\renewcommand{\i}{\mkern 1mu{\rm i}\mkern 1mu }  
\newcommand{\theorem}[1]%
{\par\addvspace{0.7\baselineskip}\par{\bf Theorem #1.}\sl}
\newcommand{\lemma}[1]%
{\par\addvspace{0.7\baselineskip}\par{\bf Lemma #1.}\sl}                                               
\newcommand{\remark}[1]{\par\addvspace{0.7\baselineskip}\par{\bf Remark #1.}}
\newcommand{\corollary}[1]%
{\par\addvspace{0.7\baselineskip}\par{\bf Corollary #1.}\sl}
\newcommand{\example}[1]{\par\addvspace{0.7\baselineskip}\par{\bf Example #1.}}
\def\VOLUME{?}
\def\YEAR{2004}
\def\ISSUE{?}

\renewcommand{\title}[1]{\thispagestyle{empty}\noindent\normalsize%
\it J. Inv.\ Ill-Posed Problems\rm, Vol.\ \VOLUME, No.\ \ISSUE,
pp.\ \pageref{FIRSTPAGE}--\pageref{LASTPAGE} (\YEAR)\par\noindent%
\copyright\ VSP \YEAR\par\vspace{6.0\baselineskip}%
\par\Large\par\noindent {\bf #1}\par\normalsize\par\label{FIRSTPAGE}}
\newcommand{\beginauthor}{\par\vspace{1.2\baselineskip}\par\noindent}
\def\endauthor{\par}
\newcommand{\rece}[1]{\par\vspace{1.2\baselineskip}\small\par\noindent%
\sl Received #1}

\newcommand{\beginabstract}{\par\vspace{1.2\baselineskip}%
\small\par\noindent{\bf Abstract} --- }
\def\endabstract{\par\normalsize\par}
\renewcommand{\div}{\mathop{\hspace{0.01pt}{\rm div}}\nolimits}   

\renewcommand{\Re}{\mathop{\hspace{0.01pt}{\rm Re}}\nolimits}
\renewcommand{\Im}{\mathop{\hspace{0.01pt}{\rm Im}}\nolimits}

\newcommand{\grad}{\mathop{\hspace{0.01pt}{\rm grad}}\nolimits}
\begin{document}
 \title{Singular value decomposition for the 2D fan-beam Radon transform
 of tensor fields}

 \beginauthor
 S.G.~Kazantsev \footnote
 {Sobolev Institute of Mathematics, Siberian Branch of Russian Academy of Sciences,
 Acad. Koptyug prosp.,~4,  Novosibirsk, 630090, Russia. E-mail: kazan@math.nsc.ru}
  and A.A.~Bukhgeim \footnote
 { Sobolev Institute of Mathematics,  Siberian Branch of Russian Academy of Sciences,
 Acad. Koptyug prosp.,~4,  Novosibirsk, 630090, Russia.
 E-mail: bukhgeim @math.nsc.ru
 \\
 \\
 \noindent  This work was partially supported by grant  RFFI\no \ 02-01-00296 and by EU grant IST-1999-29034.
}
 \endauthor

 {\rece 10.10. 2001}

\beginabstract
 In this article we study the fan-beam Radon transform ${\cal D}_m $
 of symmetrical solenoidal 2D tensor fields of arbitrary rank $m$
 in a unit disc $\mathbb D$ as the  operator, acting
 from  the object space ${\mathbf L}_{2}(\mathbb D; {\bf S}_m )$ to
 the data space  $L_2([0,2\pi)\times[0,2\pi)).$
 The orthogonal polynomial basis ${\bf s}^{(\pm m)}_{n,k}$
 of solenoidal tensor fields on the disc $\mathbb D$ was
 built with  the help of Zernike polynomials and
 then a singular value decomposition (SVD) for the
 operator ${\cal D}_m $  was obtained. The inversion formula
 for the fan-beam tensor transform ${\cal D}_m $  follows from this decomposition.
 Thus obtained inversion formula can be used as a tomographic filter for splitting
 a known tensor field into potential and solenoidal parts.
  Numerical results are presented.
\endabstract

\section{Introduction}

 \noindent The problem of determining vector or tensor field from the
 integral information arises in various applications, for instance in ultrasound probing
 of fluid or gas flows and deformed elastic media. In the first
 case it's required to determine the velocity vector field in the
 flow and in the second case --- the stress tensor field.

 One of the most complete monograph of the tensor tomography is
 \cite{Shar}. Reversibility and stability of different kinds of
 transforms of tensor fields on the Riemannian  manifolds
 are studied there. In \cite{BH}, \cite{Nor} the solution of the vector
 tomography problem is reduced to the scalar Radon problem.
 An approximate solution of the vector and tensor (of rank 2) tomography problem is given in
 \cite{BDS02}, \cite{DK02} with the help of polynomial
 non-orthogonal basis.

 More information and references about
 vector and tensor tomography problems  are given in
 \cite{Bo01}, \cite{Sch01}, \cite{St99},  \cite{Ver00}.
 \par In this article we derive an inversion formula on the basis of
 singular value decomposition (SVD) for the fan-beam transform of tensor fields.
 To this end the orthonormal polynomial  basis of solenoidal
 tensor fields, supported in unit disk,  are built from Zernike polynomials.
 In the scalar case thus obtained SVD corresponds to the known SVD for the Radon
 transform in the classical (parallel) formulation
 \cite {Corm63}, \cite {Corm64}, \cite{Herm80}, \cite {Mar74}, \cite{M90}.

 Unlike the scalar case, Radon transform of tensor fields has a
 non-zero kernel and it's possible to reconstruct uniquely (without
 additional information) only the solenoidal part of a tensor
 field, so the inversion formula can be used as a tomographic
 filter for splitting a known tensor field into potential and solenoidal parts.

 This article is organized as follows:
 In Section 2 we formulate  the problem of 2D tensor tomography.
 In Section 3 we review those part of the tensor fields theory
 that are needed in this paper.
 Section 4 contain a novel properties of Zernike polynomials.
 Sections  5 is  devoted to  the  orthogonal polynomial basis in
 the space of solenoidal (divergence free) tensor fields
 and   a singular value decomposition (SVD) for the tensor
 tomograph problem.
  A short description of the implementation issues and numerical
 tests are presented in Section 6.

 \section{Formulation of the problem}
 Let us consider the Cartesian coordinate system
 $(x^1,x^2)$ on the  plane $\mathbb R^2$ and
 let ${\bf T}_m$ denote for $m=0,1,...$ the space of all real-valued $m$-{\it covariant
 tensors}
 $${\bf a}:=a_{i_1...i_m}dx^{i_1}\otimes dx^{i_2}
 \otimes ...\otimes dx^{i_m} \ \text{or}  \
 {\bf a}=  \{ a_{i_1...i_m}, \ i_1,...,i_m =1,2 \},$$
 where $\otimes$ is the tensor product and $a_{i_1...i_m}$
 are the  components of ${\bf a}$ in the Cartesian basis $(x^1,x^2).$
 Here and throughout we imply the summation convention.
 By ${\bf S}_m$ we denote the subspace of {\it symmetric} $m$-covariant tensors and
 there exists a canonical projection $\sigma: {\bf T}_m \to {\bf S}_m$
 (called symmetrization) onto this space defined by the equation
 \begin{equation}
 (\sigma {\bf a} )_{i_1...i_m}:=\frac {1}{m!}\sum_{\pi\in \Pi_m}
 a_{i_{\pi(1)}...i_{\pi(m)}}, \label{simm}
 \end{equation} where $\Pi_m$ is the group of all permutations of degree $m.$
 A symmetric $m$-covariant tensor  ${\bf a}=\{ a_{i_1...i_m}, i_1,...,i_m =1,2\}$
 has only $m+1$ independent components which we denoted by $ a_k,$ so that
 \begin{equation}
 a_k:=a_{\underbrace{1...1}_k\underbrace{2...2}_{m-k}}, \
 (k=0,...,m). \label{neskomp}
 \end{equation}
 \par
 We will always denote vector and tensor fields and any related quantities
 such as functional spaces by boldface characters.
 \par
 Let $\mathbb D:=\{(x^1,x^2)\in {\mathbb R^2} \ \Big | \  (x^1)^2+(x^2)^2<1 \}$ be a unit
 disc on the plane ${\mathbb  R}^2.$
 The symmetric $m$-covariant  tensor field ${\bf a}(x^1,x^2)$ defined on $\mathbb D$
 can be treated as a mapping
 $${\bf a}: {\mathbb D} \rightarrow {\bf S}_m, \ {\bf a}(x^1,x^2)= \{ a_{i_1...i_m}(x^1,x^2),
  \ i_1,...,i_m =1,2 \}.$$
 \par
 {\it The fan-beam Radon transform} ${\cal D}_m$ of tensor field ${\bf a}(x^1,x^2)$
 is defined by
 \begin{equation}
  [{\cal D}_m {\bf a}](\beta, \varphi)
  := \mathop {\int}_{0}^{2\cos(\beta-\varphi) }\theta^{i_1}
  \cdot\theta^{i_2}\cdot...\cdot\theta^{i_m}
  a_{i_1...i_m}(\cos \beta-l\cos \varphi, \sin \beta-l\sin \varphi ) {\d}l,
 \label{fanT}
 \end{equation}
 where $\beta \in [0,2\pi),  \ \bm{\theta}=\begin{pmatrix} \theta^1  \\ \theta^2 \end{pmatrix}
 = \begin{pmatrix} \cos \varphi  \\  \sin \varphi \end{pmatrix},
 \  |\beta-\varphi |\leq \dfrac {\pi}{2}.$
 \\
 The difference between {\it the parallel-beam} and {\it the fan-beam
 geometry} is shown in figure 1.
 \begin{figure}[h!]
 \begin{center}
 \includegraphics[scale=0.55, bb=0 0 734 324] {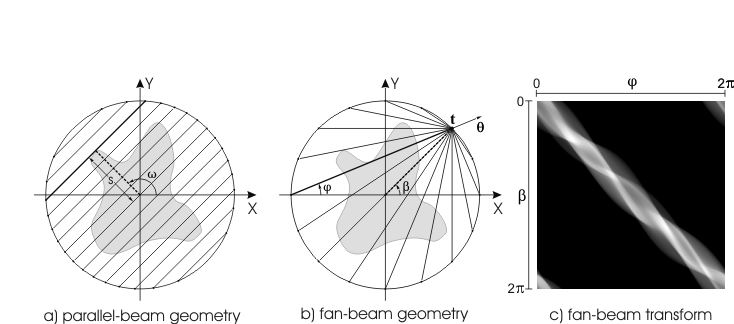}
 \caption{\footnotesize Left: parallel-beam scanning geometry. Middle: fan-beam scanning geometry.
 Right: an example of
 the fan-beam transform (the  data function or  sinogram) $f(\beta,\varphi),$ the angle $\beta$
  defines {\it the vertex point} of the fan-beam
 projection $f(\beta,\cdot)$  and the angle $\varphi$ defines
 the direction of scanning.}
 \end{center}
 \end{figure}
 \par  For  $|\beta-\varphi |> \dfrac {\pi}{2}$ we complete the definition of the
 fan-beam transform
 (\ref {fanT}) with the condition
 \begin{equation}\label{dcond}
 [{\cal D}_m{\bf a}](\beta,\varphi):=(-1)^{m+1}[{\cal D}_m{\bf a}](\beta,\varphi+\pi).
 \end{equation}
 Note, that the case $m=0$ corresponds to  the fan-beam Radon transform ${\cal D}_0\equiv {\cal D}$
 of a scalar function $a(x^1,x^2)$.
 \par Now, the problem is to recover the unknown tensor field ${\bf a}(x^1,x^2)$
 in the unit disc $\mathbb D$ from the data function  $f(\beta,\varphi),$ see Figure 1c,
 such that
 $$[{\cal D}_{m} {\bf a}](\beta,\varphi)=f(\beta,\varphi), \
 (\beta,\varphi) \in [0,2\pi) \times [0,2\pi). $$

 This problem will be solved here by the SVD-method.

 \section{Preliminaries}
 In this section, we  introduce the definition of SVD method
  and then review some facts from vector and
 tensor analysis \cite{Shar} and, in particular, consider real-valued tensor fields in
 complex coordinates (variables) \cite{V}.
 We define here some functional spaces of tensor fields ---
 ${\bf  L}_2(\mathbb D; {\bf S}_m),$ for example, and
 also establish the notations that will be
 used in the sequel.

 \subsection{ Singular value decomposition (SVD)}%
Now we define the concept of a {\it singular value decomposition},
 see \cite{Nat86}, \cite{NW01}, \cite{M90}, \cite{Qui83}.
  Let  $U$ and $V$ be Hilbert spaces,
  and  $A$ be a compact linear operator
 from  $U$ to  $V, \ A\in {{\cal L}(U,V)}. $ Then there exists a sequence $\{ \sigma_k \}_{k \geq 1}$
 of positive numbers, monotonically tending to zero (or a finite sequence)
 and two orthonormal systems $\{u_k \}_{k \geq 1}\subset U,\
 \{v_k\}_{k\geq 1}\subset V,$
 such that for all $u\in U$ we have a singular value decomposition
 $$Au=\sum_{k=1}^{\infty}\sigma_k (u,u_k)_Uv_k, \ Au_k=\sigma_k v_k, \
 \sigma_1\geq \sigma_2\geq ...>0.$$
 The adjoint of $A$ is given by
 $$A^*v=\sum_{k=1}^{\infty}\sigma_k (v,v_k)_Vu_k, \ A^*v_k=\sigma_k u_k$$
 and the generalized inverse of $A$ is
 $$A^+v=\sum_{k=1}^{\infty}\sigma_k^{-1} (v,v_k)_Vu_k.$$
 Operator  $A^+$ can be unbounded, so one can use a truncated SVD for its regularization
 $$T_\gamma v=\sum _{k \leq 1/\gamma}\sigma_k ^{-1}(v,v_k)_Vu_k,$$
 where $\gamma$ is the parameter of regularization.
 SVD is one of the methods for solving ill-posed problems and
 it allows to characterize the range of the operator, invert it and estimate an incorrectness
 of the corresponding inverse problem.
 \subsection{Tensor fields in complex coordinates}
 Let's identify $\mathbb R^2$ with the complex plane $\mathbb C$
 by the usual way
 $$z^1\equiv z:=x^1+{\i}x^2, \ z^2 \equiv
 \overline{z}:=x^1-{\i}x^2, \ {\i}^2=-1.$$
 Let ${\bf a}=\{ a_{i_1...i_m}(x^1,x^2) \}$ be an $m$-covariant real-valued
 tensor field in Cartesian coordinates $(x^1, x^2),$  then in complex
 coordinates or variables $(z,\overline z)$
 it will have new components $A_{i_1...i_m}(z,\overline{z}),$
 which are formally expressed by {\it the covariant tensor law}
 \begin{equation}
 A_{i_1...i_m}(z,\overline{z})=
 \frac {\partial x^{s_1}}
 {\partial z^{i_1}}...
 \frac {\partial x^{s_m}}{\partial z^{i_m}}a_{s_1...s_m}(x^1,x^2),
 \label{compkoor}
 \end{equation}
 where the Jacobian matrix is
 \begin{equation}
 J\equiv(J^i_j):= \begin{pmatrix}
  \dfrac {\partial z^1}{\partial x^1} & \dfrac {\partial z^1}{\partial x^2} \\
 \\
 \dfrac {\partial z^2}{\partial x^1}  & \dfrac {\partial z^2}{\partial x^2} \\
   \end{pmatrix}= \
   \begin{pmatrix}
  1 & \ {\i} \\
  1 & - {\i} \\
  \end{pmatrix} \nonumber
 \end{equation} \\
 and the inverse matrix  of it is
 \begin{equation}
  J^{-1}=
 \begin{pmatrix}
 \dfrac {\partial x^1}{\partial z^1} &\dfrac {\partial x^1}{\partial z^2} \\
 \\
 \dfrac {\partial x^2}{\partial z^1} &\dfrac {\partial x^2}{\partial z^2} \\
 \end{pmatrix}= \
 \frac {1}{2} \begin{pmatrix}
 \ 1    &  1   \\
 - {\i} & {\i} \\
 \end{pmatrix}. \nonumber
 \end{equation}

 Here the  formal partial derivatives with respect to
 $z^1$ and $z^2$ are defined in the usual way
 \begin{eqnarray}
 \frac {\partial }{
 \partial z^1}\equiv \frac {\partial }{ \partial z}:= \frac {1}{2}\left
 ( \frac
 {\partial }{ \partial x^1}-{\i}\frac {\partial }{ \partial x^2} \right
 ), \
 \frac {\partial }{\partial z^2}\equiv
 \frac {\partial }{ \partial \overline z}:=
 \frac {1}{2}\left ( \frac {\partial }{ \partial x^1}+{\i}\frac
 {\partial }{ \partial x^2} \right ).
 \label{part}
 \end{eqnarray}

 We shall write transformation   (\ref{compkoor}) as
 $${\bf a}=\{ a_{i_1...i_m}(x^1,x^2) \} \rightarrowtail {\bf A}=\{A_{i_1...i_m}(z,\overline{z}) \}.$$
 From now on small letters will be used to denote tensor fields in
 the initial Cartesian coordinate system $(x^1,x^2)$ and capital
 letters will be used for the same tensor fields in complex
 coordinates $(z,\overline{z}).$

 An inverse relationship also takes place
 \begin{equation}
 a_{i_1...i_m}(x^1,x^2)=\dfrac {\partial z^{s_1}}{\partial x^{i_1}}...
 \dfrac {\partial z^{s_m}}{\partial x^{i_m}}A_{s_1...s_m}(z,\overline{z}),
 \label{rtocom}
 \end{equation}
 and we shall also write this as
 $${\bf A}=\{A_{i_1...i_m}(z,\overline{z}) \} \rightarrowtail
 {\bf a}=\{ a_{i_1...i_m}(x^1,x^2) \}.$$

 A symmetric $m$-covariant  tensor ${\bf A}$  could also be given by its
 components $A_k$
 \begin{equation}
 A_k:=A_{\underbrace{1...1}_k\underbrace{2...2}_{m-k}}, \
 (k=0,...,m)
 \label{Ak}
 \end{equation}
 and subject to the conditions
 \begin{equation}
 A_k=\overline{A}_{m-k}, \ (k=0,...,m).
 \label{cond}
 \end{equation}
 So we may image the symmetric tensor as pseudovector,
 expanding the one as a column array for convenience, that
 the following notations will be used
 \begin{equation}\label{pseudo}
 {\bf a}=\begin{pmatrix}  a_m  \  \\
 a_{m-1}\\ ...\\  a_1 \ \\  a_0 \
 \end{pmatrix}, \
 {\bf A}=\begin{pmatrix}  A_m \  \\
 A_{m-1}\\ ...\\ A_1 \ \\ A_0 \
 \end{pmatrix}
 \end{equation}

 Taking into account the tensor law (\ref{compkoor}), (\ref{rtocom})
 we get the formulae that link independent components (\ref{neskomp}) and
 (\ref{Ak}) in  pseudovectors (\ref{pseudo})
 \begin{eqnarray} \label{a_k}
 a_k&=&(-{\i})^{m-k}\sum_{p=0}^{m-k} \sum_{q=0}^{k}
 C^{p}_{m-k}C^q_{k}(-1)^p A_{p+q}, \\
 \label{A_k}
 A_k&=&
 \frac {{\i}^{m-k}}{2^m}\sum_{r=0}^{m-k} \sum_{s=0}^{k}
 C^r_{m-k}C^s_{k}(-{\i})^{k+r-s} a_{r+s},
 \end{eqnarray}
 where $k=0,1,...,m$ and $C^i_j$ are binomial coefficients.

\subsection{Metric tensor ${\bf G} $ and pointwise inner product in complex coordinates}

 On parity with covariant components of the tensor we shall also use its
 {\it contravariant components}. In Cartesian coordinates $(x^1,x^2)$
 covariant and contravariant components $g_{ij}$ and $g^{ij}$
 of the metric tensor $\bf g$ are the same
 \begin{equation}
 {\bf g}=\{g_{ij} \}= \{g^{ij} \}=
 \left \{
 \begin{array}{cc}
  1 & 0 \\
  0 & 1
 \end{array}
 \right \}.  \label{gij}
 \end{equation}

 Thus contravariant components of the tensor ${\bf a}$ coincide with
 its corresponding covariant components,  $a_{i_1...i_m}=a^{i_1...i_m}.$
 The  pointwise inner product $\langle \cdot, \cdot \rangle$
 on ${\bf S}_m$
 induced by the Euclidian metric
 ${\bf g}$  (\ref{gij})    is defined by the formula

 \begin{equation}\label{scal}
 \langle {\bf a},~{\bf b} \rangle:=a^{i_1i_2...i_m }b_{i_1i_2...i_m }.
 \nonumber
 \end{equation}
 In complex coordinates $(z,\overline z)$ the metric tensor ${\bf G}$
 has the following covariant
 \begin{equation}
 \{G_{ij}\}=  \left \{ \begin{array}{cc}
 0 & 1/2 \\
 1/2 & 0
 \end{array} \right \}
 \label {Gcov} \nonumber
 \end{equation}
 and contravariant components
 \begin{equation}
 \{G^{ij} \}=\left \{
 \begin{array}{cc}
 0 & 2 \\
 2 & 0
 \end{array} \right \}.
 \label{Gcontr}
 \end{equation}
 Contravariant components of tensor ${\bf A}$ in complex coordinates
 are obtained by raising indexes with contravariant components
 of the metric tensor (\ref{Gcontr})
 \begin{equation}
 A^{i_1i_2...i_m } =G^{i_1j_1}G^{i_2j_2}...G^{i_mj_m} A_{j_1j_2...j_m}
 \nonumber
 \end{equation}
 and the pointwise inner product of tensor fields is evaluated by formula
 \begin{eqnarray}
 \langle {\bf A},{\bf B} \rangle= A^{i_1i_2...i_m }B_{i_1i_2...i_m }=
 A_{i_1i_2...i_m}B^{i_1i_2...i_m}. \label{Scal}
 \nonumber
 \end{eqnarray}

 If tensors ${\bf A}=\{A_k \}$ and ${\bf B}=\{ B_k \}$ are considered as
  pseudovectors
 in complex coordinates
 then their pointwise
 inner product will be equal to
 \begin{equation}\label{pscal}
 \langle {\bf A},~{\bf B}\rangle=2^m\sum_{k=0}^{m}C^k_mA_kB_{m-k}.
 \end{equation}
 The pointwise norm of tensor ${\bf A}$ then will be
 \begin{equation}
 |{\bf A}|^2=2^m\sum_{k=0}^{m}C^k_m|A_k|^2.
 \label{norm}
 \end{equation}
 It is clear that the pointwise inner product is invariant, i.e.
 if ${\bf a}  \rightarrowtail {\bf A}$ and ${\bf b}  \rightarrowtail {\bf B},$ than
 \begin{equation}
 \langle {\bf a},~{\bf b}\rangle=\langle {\bf A},~{\bf B}\rangle.
 \label{invsc}
 \end{equation}
\subsection{The space of integrable tensor fields
 ${\bf L}_2(\mathbb D; {\bf S}_m)$}
 Let ${\mathbf L}_2(\mathbb D; {\bf S}_m)$ denote a Hilbert space
 comprising real-valued symmetric
 $m$-covariant tensor fields on $\mathbb D$
 with the inner product, denoted by $\langle \langle \cdot,\cdot
 \rangle \rangle $
 \begin{equation}\langle \langle {\bf a},~{\bf b}
  \rangle \rangle \equiv
  \langle \langle {\bf  a},~{\bf b}
  \rangle \rangle_{{\mathbf  L}_2(\mathbb D; {\bf S}_m)}:=
  \iint\limits_{\mathbb D} \langle {\bf  a}(x^1,x^2),~{\bf b}(x^1,x^2)\rangle  {\d}V^2,
 \ {\d}V^2={\d}x^1 \wedge  {\d}x^2
 \nonumber
 \end{equation}
 and the  finite norm $||\cdot||$
 \begin{equation} ||{\bf a}||^2\equiv
 ||{\bf a}||^2_{{\mathbf  L}_2(\mathbb D; {\bf S}_m)}:=
 \iint\limits_{\mathbb D} \langle {\bf a}(x^1,x^2),~{\bf a}(x^1,x^2)\rangle
 {\d}V^2.
 \nonumber
 \end{equation}
 In complex coordinates for ${\bf a} \rightarrowtail {\bf A}$ and  ${\bf b} \rightarrowtail {\bf B}$ we have
 $$\langle \langle {\bf A},~{\bf B}\rangle \rangle=
 \iint\limits_{\mathbb D} \langle {\bf A}(z,\overline z),~{\bf B}(z,\overline z)\rangle {\d}V^2,
 \ {\d}V^2=\frac {{\d}z\wedge {\d}\overline z}{-2{\i}}.$$
 By virtue of invariance of inner product (\ref{invsc}) the
  following equalities take place
   $$\langle \langle {\bf a},~{\bf b}\rangle \rangle
 =\langle \langle {\bf A},~{\bf B}\rangle \rangle, \ ||{\bf a}||=||{\bf A}||. $$

\subsection{Differential operations on symmetric tensor fields}
 We shall denote the class of real-valued $m$-covariant symmetric tensor fields
 ${\bf a}=\{a_{i_1...i_m}(x^1,x^2)\},$
 whose all components
 are functions from $C^k( \mathbb D),\ 1 \leq k \leq \infty$ by
 ${\mathbf C}^k(\mathbb D ; {\bf S}_m).$
 A subset of ${\mathbf C}^k( \mathbb D; {\bf S}_m)$
 whose finite support is contained in $\mathbb
 D$ will be denoted by ${\mathbf C}_0^{k}(\mathbb D ; {\bf S}_m).$

 The operator of {\it covariant differentiation}  $\nabla$
  (in the vectorial case $ \equiv \grad $)
 $$\nabla:{\mathbf C}^{\infty}(\mathbb D; {\bf S}_{m})
 \to {\mathbf C}^{\infty}(\mathbb D; {\bf S}_{m+1})$$
 in Cartesian coordinate system  $(x^1,x^2)$ is defined by equation
 $$\nabla {\bf a}:= \{ a_{i_1...i_m;j} \}=
 \left  \{ \frac {\partial
 a_{i_1...i_m}} {\partial x^j}, \ j=1,2   \right  \}.
 $$
 The  covariant differentiation  $\nabla$  operates on any tensor field of rank $m \geq 0$
 and produces a tensor field that
 is one rank higher. For example,
 the gradient of a (co)vector field  is a second rank tensor field.

 In complex coordinates we have
  $$\nabla {\bf A}
  =\{ A_{i_1...i_m;j} \}=\left  \{ \frac {\partial A_{i_1...i_m}} {\partial z^j}, \ j=1,2\right  \},$$
 where formal partial derivatives with respect to
 $z^1$ and $z^2$ are defined
 by (\ref{part}).
 \par The operator of {\it divergence} $\delta$ (in the vectorial case $ \equiv \div$)
 $$\delta :{\mathbf C}^{\infty}(\mathbb D; {\bf S}_m) \to
 {\mathbf C}^{\infty}(\mathbb D ; {\bf S}_{m-1} )
 $$
 in Cartesian coordinate system $(x^1,x^2)$ is defined by
 $$\delta {\bf a}:= \{
 a^{i_1i_2...i_{m-1}s}_{;x^{j}} \}=
 \left  \{ \frac {\partial a^{i_1i_2...i_{m-1}j}}{\partial x^{j} }, \ j=1,2  \right \}.
 $$
 The divergence $\delta$ can operate on any tensor field of rank $m\geq 1$
 and above produces a tensor that is one rank lower. For example,
 the divergence of a second
 rank tensor field  is a (co)vector field.
 \par In complex variables the divergence is calculated with the help of
  contravariant components $G^{ij}$ of the metric tensor (\ref{Gcontr})
 \begin{eqnarray}
 \delta {\bf A}&=& \left   \{ A_{i_1i_2...i_m;z^s} G^{i_ms} \right \} \nonumber \\
 &=&
 \left \{
 2\frac {\partial A_{i_1i_2...i_{m-1}2}}{\partial z}+
 2\frac {\partial A_{i_1i_2...i_{m-1}1}}{\partial \overline{z}}, \   i_1,...,i_{m-1}=1,2
 \right \}. \label{cdiv}
 \end{eqnarray}
 A smooth tensor field  ${\bf a} \in {\mathbf C}^k(\mathbb D; {\bf S}_m)$
 is called {\it solenoidal} if its divergence equals to zero.
 The condition for the tensor field ${\bf A}$ to be solenoidal
 can be expressed in complex
 coordinates in terms of its independent components
 $A_0,...,A_m$
 \begin{eqnarray}
 \begin{tabular} {l} $ \left \{\begin{array} {ll}
 (A_{0})_z+(A_{1})_{\overline{z}}&=0  \\
 ...&   \\
 (A_{k})_z+(A_{k+1})_{\overline{z}}&=0 \\
 ...&  \\
 (A_{m-1})_z+(A_{m})_{\overline{z}}&=0
 \end{array} \right.$
 \end{tabular}
 \ \text {or} \ \
 \frac {\partial }{ \partial \overline z}
 \begin{pmatrix} A_{m}  \\ A_{m-1}\\ ...\\ A_2\\ A_1
 \end{pmatrix}+\frac {\partial}{ \partial z}
 \begin{pmatrix} A_{m-1}  \\ A_{m-2}\\ ...\\ A_1\\ A_0
 \end{pmatrix}=0.
 \label{sol}
 \end{eqnarray}
 \par The next differential operation on the
 symmetric tensor fields is  {\it the symmetric inner differentiation}  $d$
 $$d : {\mathbf C}^{\infty}(\mathbb D;
 {\bf S}_{m-1} )\to {\mathbf C}^{\infty}( \mathbb D ; {\bf S}_{m}),
 $$
 defined in the following way
 $$d:=\sigma \nabla,
 $$
 where $\sigma$ is the symmetrization operator (\ref{simm}).
 \par A tensor field  ${\bf a}\in {\mathbf C}^{\infty}
 (\mathbb D ; {\bf S}_m)$
 is called a smooth  {\it potential field}, if for
 some tensor field  ${\bf v}\in {\mathbf C}_0^{\infty}
 (\mathbb D ; {\bf S}_{m-1})$ with
 boundary condition ${\bf v}|_{\partial \mathbb D }=0 $
 we have ${\bf a}=d{\bf v}$ and ${\bf v}$ is the  potential.
 \par
 The symmetric inner differentiation $d$ in complex variables
 is calculated in the following manner. If $ {\bf a} =d{\bf v}$ and $ {\bf a} \rightarrowtail{\bf A}, \
 {\bf v} \rightarrowtail{\bf V}$ then ${\bf A}=d{\bf V}$ and
 \begin{equation} \label{poten}
 A_k=\frac {m-k}{m}\frac {\partial V_{k}}{ \partial \overline z}
 +\frac {k}{m}\frac {\partial  V_{k-1} }{ \partial  z}, \   (m \geq 1, \ k=0,1,...,m) .
 \end{equation}
\subsection{Orthogonal decomposition of the space
 ${\bf L}_{2}(\mathbb D; {\bf S}_m)$ into the sum of solenoidal and potential parts}
 Operators $d$ and $-\delta$ are formally conjugate and for
 a bounded region ${\mathbb G}$ with a piecewise-smooth boundary $\partial {\mathbb G}$
 the Gauss-Ostrogradsky formula takes place
 \begin{eqnarray}
 \iint\limits_{\mathbb G}  (\langle d {\bf v},~{\bf a} \rangle+\langle
 {\bf v},\delta {\bf a}\rangle ){\d}V^2
 =\mathop{\int}_{\partial {\mathbb G}} \langle i_{\bm{\nu}} {\bf v},~{\bf a}\rangle {\d}V^1,
 \label{Gauss}
 \end{eqnarray}
 where
 ${\bf a}\in {\bf S}_m, \ {\bf v}\in {\bf S}_{m-1}$ are smooth tensor
 fields, and $\bm{\nu}= \{ \nu_1,\nu_2 \} \in {\bf S}_1$ is a unit covector of outward normal
 to the boundary $\partial {\mathbb G},$ and $i_{\bm {\nu}}$ is the operator of
 {\it symmetric multiplication} with the covector $\bm {\nu}$
 $$i_{\bm {\nu}} : {\bf S}_m \to {\bf S}_{m+1},$$
 which is defined by the equation
 $$(i_{\bm {\nu}} v)_{i_1...i_mi_{m+1}}
 :=\sigma(\nu_{i_1} v_{i_2...i_{m+1}}).$$
 In terms of the Gauss-Ostrogradsky formula (\ref {Gauss})
 we can define that a tensor field
 ${\bf a}\in {\mathbf L}_2({\mathbb D}; {\bf S}_m)$ is solenoidal if the following
 equation takes place
 \begin{eqnarray}
 \langle \langle d{\bf v},~{\bf a} \rangle \rangle=
 \iint\limits_{\mathbb D} \langle d{\bf v},~{\bf a}\rangle {\d}V^2=0
 \label{defsol}
 \end{eqnarray}
 for all smooth tensor fields
 ${\bf v}(x^1,x^2) \in {\mathbf C}_0^{\infty}(\mathbb D; {\bf S}_{m-1}).$
 \par  We denote by ${\mathbf H}({\mathbb D}; {\bf S}_m, \delta)$
 the graph space of $\delta $ over ${\mathbf L}_2({\mathbb D}; {\bf S}_m),$ i.e.
 $${\mathbf H}({\mathbb D}; {\bf S}_m, \delta):=\{
  {\bf u}\in {\mathbf L}_2({\mathbb D}; {\bf S}_m) \ \Big| \ \delta {\bf u}
  \in{\mathbf L}_2({\mathbb D}; {\bf S}_m) \}.
  $$
  It is a Hilbert space under the graph norm
  $$\langle \langle {\bf u}, {\bf v}\rangle \rangle_
   {{\mathbf H}({\mathbb D}; {\bf S}_m, \delta)}
  :=\langle \langle {\bf u}, {\bf v}\rangle \rangle +
  \langle \langle \delta{\bf v}, \delta{\bf v}\rangle \rangle, \
 ||{\bf u}||^2_
   {{\mathbf H}({\mathbb D}; {\bf S}_m, \delta)}
  :=|| {\bf u}||^2 +
   ||\delta{\bf u}||^2.
  $$
  Finally, we define subspace of solenoidal tensor fields
  (i.e. which  satisfies the equation (\ref{defsol}))
 $${\mathbf H}({\mathbb D}; {\bf S}_m, \delta=0):=\{
  {\bf a}\in{\mathbf H}({\mathbb D}; {\bf S}_m, \delta)\Big| \ \delta {\bf a}=0
  \}
  $$
 and it is clear that this subspace  is a completion  of the set of smooth solenoidal tensor fields
 with respect to the norm $|| \cdot||$ of  ${\mathbf L}_{2}(\mathbb D ;  {\bf S}_m ).$

By ${\bf H}_N(\mathbb D ;{\bf S}_m,\delta=0) $
 we denote the finite-dimensional  subspace of polynomial (of degree at most $N$)
 solenoidal $m$-covariant tensors fields.
 Then we have
 $${\bf H}_0(\mathbb D ;{\bf S}_m,\delta=0) \subset{\bf H}_1(\mathbb D ;_m,\delta=0)
  \subset...\subset
 {\bf H}_N({\mathbb D ;\bf  S}_m,\delta=0) \subset ...\subset L_2(\mathbb D; {\bf S}_m)$$
 and
 $${\bf H}(\mathbb D ;{\bf S}_m,\delta=0)={\bf clos} \left ( \bigcup_{N=0}^{\infty }
 {\bf H}_N({\mathbb D ;\bf S}_m,\delta=0) \right ), $$
 where {\bf clos}
  means the closure in
 ${\mathbf L}_{2}(\mathbb D ;  {\bf S}_m ).$

 \par
 It is well known, see  \cite{ DL}, \cite {W}, that a vector field can be
 represented as a sum of solenoidal and potential vector fields.
 The classical result in this direction belongs to H. Weyl and is
 connected with the decomposition of the
 $L_2$ space of vector fields into the orthogonal sum of solenoidal and potential fields.
 The analogous result is true for tensor fields, see \cite{Ca81},
 \cite{G-MA94}, \cite{Shar}. Namely, for   ${\bf u}\in {\mathbf L}_{2}(\mathbb D ;
 {\bf S}_m )$ we have
 $$ {\bf u}={\bf a}+d{\bf v}, \ \langle \langle {\bf a},~{\bf d{\bf v}}=0  \rangle \rangle,
 $$
 where
  ${\bf a}$
  is a solenoidal tensor field and
   ${\bf v}\in {\mathbf H}^1_0( \mathbb D ; {\bf S}_{m-1}).$
Or, in another words the orthogonal decomposition
 $${\mathbf L}_{2}(\mathbb D ;
 {\bf S}_m )={\bf H}({{\mathbb D};\bf S}_m,\delta=0)  \oplus d
 {\mathbf H}^1_0( \mathbb D ; {\bf S}_{m-1})
 $$
 takes place, where the Sobolev space ${\mathbf H}^1_0( \mathbb D ; {\bf S}_{m-1})$
 is a completion  of the space of smooth  tensor fields ${\mathbf C}_0^{1}(\mathbb D ; {\bf S}_{m-1})$
 with respect to the Sobolev norm $|| \cdot||_1,$
 corresponding to the scalar product $\langle \langle {\cdot},~{\cdot }\rangle \rangle_1$
 that is defined by the formula
 $$\langle \langle {\bf u},~{\bf v }\rangle \rangle_1=\langle \langle {\bf u},~{\bf v }\rangle \rangle+
 \langle \langle {\nabla \bf u},~{ \nabla\bf v }\rangle \rangle.
 $$

\subsection{Fan-beam Radon transform ${\cal D}_m$ of tensor fields in complex variables}%

 Let's assume that some (constant) vector field is given in Cartesian coordinates
 $$\bm {\theta}=\begin{pmatrix} \theta^1  \\ \theta^2 \end{pmatrix}
 = \begin{pmatrix} \cos \varphi  \\  \sin \varphi \end{pmatrix},$$
 then according to the tensor law for contravariant components its
 representation in complex coordinates
 will look like
 \begin{eqnarray}
 \bm {\theta} \rightarrowtail {\mathbf \Theta}, \
 \Theta^j=\frac {\partial z^j}{\partial x^s}\theta^s, \
 {\bf \Theta}=\begin{pmatrix} \Theta^1  \\ \Theta^2
 \end{pmatrix}
 = \begin{pmatrix}e^{{\i}\varphi} \\e^{-{\i}\varphi} \end{pmatrix}.
 \label{Psij}
 \nonumber
 \end{eqnarray}
 Then we denote by $\bm {theta}^m$ the tensor product
 $$\bm{\ theta}^m:=\underbrace{\bm {\theta}\otimes  \bm {\theta}\otimes ...
 \otimes \bm{\theta}}_m
 = \{ \theta^{j_1}\cdot \theta^{j_2}\cdot...\cdot\theta^{j_m} \},$$
 and  $\bm{\theta}^m$ will be an $m$-contravariant tensor
 in Cartesian coordinates and in complex coordinates we have
 the tensor product
 $${\bf \Theta}^m:=\underbrace{{\bf \Theta} \otimes{\bf \Theta} \otimes...\otimes{\bf \Theta}}_m=
 \{ \Theta^{j_1}\cdot\Theta^{j_2}\cdot...\cdot\Theta^{j_m} \}.
 $$
 It is clear that $\bm {\theta}^m \rightarrowtail {\bf \Theta}^m. $
 As soon as the inner product of tensors is invariant
 (\ref{invsc}), we get
 $$\langle {\bf a},~\bm{\theta}^m\rangle =a_{j_1...j_m}
 \theta^{j_1}\cdot\theta^{j_2}\cdot...\cdot\theta^{j_m}=
 \langle {\bf A},~{\bf \Theta}^m\rangle =A_{j_1...j_m}
 \Theta^{j_1}\cdot \Theta^{j_1}\cdot...\cdot\Theta^{j_m}. $$
 Thus we can evaluate the fan-beam transform
 (\ref{fanT})  through the components of the tensor ${\bf A}(z,\overline z)$
 \begin{eqnarray}
 [{\cal D}_m{\bf a}](\beta,\varphi)
 = \mathop {\int}_{0}^{2\cos(\beta-\varphi) }
 \langle\bm{\theta}^m, {\bf a}(\cos \beta-l\cos \varphi, \sin \beta-l\sin \varphi )\rangle {\d}l
 \nonumber \\
 = \mathop{\int}\limits_{\tau(t,\varphi)}^{\quad t}\langle{\bf \Theta}^m,
 ~{\bf A}(\zeta,\overline \zeta)\rangle|{\d}\zeta|=
 \mathop{\int}\limits_{\tau(t,\varphi)}^{\quad t}
 \Theta^{j_1}...\Theta^{j_m}
 A_{j_1...j_m}(\zeta,\overline{\zeta})
 |{\d}\zeta|, \label{comfanT}
 \end{eqnarray}
  where  $t=e^{{\i}\beta}, \ \tau (t,\varphi)=-\overline t e^{2{\i}\varphi}=-e^{{\i}(2\varphi -\beta)}, \
  \beta \in [0,2\pi), \
 \varphi \in \left [\beta -\dfrac {\pi}{2}, \beta+ \dfrac {\pi}{2} \right ].$
 \par
 Here and in the sequel, we use notation
  $$\mathop\int \limits_{z_1}^{z_2}... |{\d}\zeta|$$ for
  a line integral along
 the line segment with end points $z_1, \ z_2 \in \overline {\mathbb D}.$
 \par
 At last we can get the fan-beam transform (\ref{comfanT}) in terms of components $A_k$
 and for $\varphi \in \left [\beta -\dfrac {\pi}{2}, \beta+ \dfrac {\pi}{2} \right ]$
 we have

 \begin{eqnarray}\label{Vtr0}
 && [{\cal D}_m{\bf a}](\beta,\varphi)=\mathop {\int}_{\tau(t,\varphi)}^t
 \sum_{k=0}^{m} C^k_me^{{\i}k\varphi}e^{-(m-k)\varphi}A_k(\zeta,\overline{\zeta})
 |{\d}\zeta|   \\
 \label{Vtr1}
 &&=\sum_{k=0}^{m}
 C^k_me^{{\i}(2k-m)\varphi}\mathop {\int}_{\tau(t,\varphi)}^t
 A_k(\zeta,\overline{\zeta}) |{\d}\zeta|=\sum_{k=0}^{m}C^k_me^{{\i}(2k-m)\varphi}[{\cal D}A_k]  \\
 &&= \sum_{k=0}^{m}C^k_me^{{\i}(m-2k)\varphi}[{\cal D}A_{m-k}]. \label{Vtr2}
 \end{eqnarray}
 \par
 Recall that for $|\beta-\varphi |\geq \frac {\pi}{2}$
 the fan-beam transform ${\cal D}_m{\bf a}$ is defined by condition (\ref{dcond}).

 Now we verify that the potential part of a tensor field is ``invisible'' for tensor transform ${\cal
 D}_m.$  Let ${\bf a}=d{ \bf v}$ and ${\bf a} \rightarrowtail {\bf A}=d{ \bf V}, \ {\bf v} \rightarrowtail {\bf V}.$
 Substituting the potential tensor (\ref{poten}) in the expansion  (\ref{Vtr0})  and
 making evident evaluations, we get
 \begin{eqnarray}
 &&[{\cal D}_m{\bf a}](\beta,\varphi)=\mathop {\int}_{\tau(t,\varphi)}^t
 \sum_{k=0}^{m} C^k_me^{{\i}k\varphi}e^{-(m-k)\varphi}
 \left (\frac {m-k}{m}\frac {\partial V_{k}}{ \partial \overline z}
 +\frac {k}{m}\frac {\partial  V_{k-1} }{ \partial z}
 \right) d |z|  \nonumber \\
 &&=\mathop {\int}_{\tau(t,\varphi)}^t
 \sum_{k=0}^{m}e^{{\i}(2k-m)\varphi}\left (
 C^k_{m-1}\frac {\partial V_{k}}{  \partial  \overline z}
 +C^{k-1}_{m-1}\frac {\partial  V_{k-1} }{ \partial z}
 \right)d|z|  \nonumber \\
 &&=\mathop {\int}_{\tau(t,\varphi)}^t \left (
 \sum_{k=0}^{m-1}e^{{\i}(2k-m)\varphi}C^k_{m-1}\frac {\partial V_{k}}{  \partial \overline z}
 +\sum_{k=1}^{m}e^{{\i}(2k-m)\varphi}C^{k-1}_{m-1}\frac {\partial  V_{k-1} }{ \partial  z}
 \right )d|z|
  \nonumber \\
&&=\mathop {\int}_{\tau(t,\varphi)}^t \left (
 \sum_{k=0}^{m-1}e^{{\i}(2k-m)\varphi}C^k_{m-1}\frac {\partial V_{k}}{  \partial  \overline z}
 +\sum_{k=0}^{m-1}e^{{\i}(2k-m+2)\varphi}C^{k}_{m-1}\frac {\partial  V_{k} }{ \partial  z}
 \right )
 d|z| \nonumber \\
&&= \sum_{k=0}^{m-1}e^{{\i}(2k-m+1)\varphi}C^k_{m-1}\mathop
 {\int}_{\tau(t,\varphi)}^t
 \label{deriv}\frac {\partial V_{k} }
 {\partial {\bf \Theta}}
  d|z| \nonumber \\
 &&= \sum_{k=0}^{m-1}e^{{\i}(2k-m+1)\varphi}C^k_{m-1}
 \left ( V_{k}(t,\overline t)- V_{k}(\tau, \overline \tau)
 \right )=0.  \nonumber
 \end{eqnarray}
  Here
  \begin{eqnarray}
  \label{derivteta}\frac {\partial }
 {\partial {\bf \Theta}}:=e^{{\i}\varphi}
 \frac{\partial}{\partial z}+e^{-{\i}\varphi}\frac{\partial}
 {\partial \overline{z}}
  \end{eqnarray}
  is the derivative  in the direction
 ${\bf \Theta}=\begin{pmatrix} \Theta^1  \\ \Theta^2 \end{pmatrix}
 = \begin{pmatrix}e^{{\i}\varphi} \\e^{-{\i}\varphi} \end{pmatrix},$
 written in the complex form
 and after integration  we take into account that the potential $\bf V$
 vanishes  on the boundary of the disc $\mathbb D$.
 \par At the end of this section  we resume, that  we consider the operator
 ${\cal D}_m$ as follows
 $${\cal D}_m:{\mathbf L}_{2}(\mathbb D; {\bf S}_m )) \to L_2([0,2\pi)\times[0,2\pi),
 $$
 and $\ker {\cal D}_m$ coincides with the space of potential fields
 $d{\mathbf H}^1_0( \mathbb D ; {\bf S}_{m-1}),$ so one can say that potential
 fields are ``invisible'' for the tensorial Radon transform ${\cal D}_m.$

 \section{Zernike polynomials}
 We will identify complex plane $\mathbb C$
 with ${\mathbb R^2}$ as above.
 Let $\mathbb D=\{ z~:~|z| <1 \}$ be the  open unit
 disc in  $\mathbb C$ and  $L_2(\mathbb D)$ denote a Hilbert space
 comprising square integrable (complex-valued) functions on
 $\mathbb D $
 with the inner product denoted by $ \langle\langle\cdot,\cdot\rangle\rangle $
 $$\langle\langle a,b\rangle\rangle
  := \iint\limits_{\mathbb D}a(x^1,x^2)\overline {b(x^1,x^2)}dV^2,
 \ {\d}V^2={\d}x^1 \wedge  {\d}x^2
 $$
 and the  finite norm $||\cdot||$
 $$
 ||a||^2:=\iint\limits_{\mathbb D} |a(x^1,x^2)|^2 dV^2.
 $$
 It is known that Zernike polynomials \cite {BW83} form a complete
 orthogonal system (basis) over the Hilbert space
 $L_2(\mathbb D).$  Let  $z=re^{{\i}\psi}\in \mathbb D, \ r=|z|, \ \psi=\arg(z).$
 Traditionally, Zernike polynomials, see \cite{BW83}, \cite{Mar74}, \cite{PR89},
 are defined by
 \begin{equation}\label{zernike}
  V_{n,l}(r,\psi ):=e^{{\i}l\psi} R_{n,|l|}(r),
 \ (-n\leq l\leq n, \ ),
 \end{equation}
 where
 \begin{eqnarray} \label{radialZ}
 R_{n,m}(r):=\sum_{p=0}^{(n-m)/2}(-1)^p
 \frac {(n-p)!}{p!\left ( \frac {n+m}{2}-p \right )
                  \left ( \frac {n-m}{2}+p \right )}r^{n-2p} \nonumber
 \end{eqnarray}
 are the so-called real-valued Zernike radial polynomials \cite{PR89},
 defined for integers $n$ and $m$ so that $0\leq m\leq n$ and
 $n-m$ is even.
 The family $R_{n,m}(r)$
 is related
 to the   Jacobi polynomials
 \begin{eqnarray} \label{RJac}
 R_{n,m}(r)=r^{m}P_{\frac {n-m}{2}}^{(0,m)}(2r^2-1),
 \nonumber
 \end{eqnarray}
 where  the  Jacobi polynomials \cite{GR80}
 are given through the Rodriguez formula
 \begin{equation} \label{jacobi}
  P_k^{(a,b)}(s):= \frac {(-1)^k}{k!2^k}(1-s)^{-a}
 (1+s)^{-b}\frac {d^k}{d s^k}\left [
 (1-s)^{(k+a)}
 (1+s)^{(k+b)}
 \right ]
 \nonumber
  \end{equation}
 for $a,b >-1; \ n=0,1,2,..., $ and $s\in {\mathbb R}$
 form a complete orthogonal
 system in the Hilbert space $L_2[-1,1]$ of square integrable
 functions on $[-1,1].$

  In this paper we will use another numbering of Zernike polynomials
  (\ref{zernike})
  and treat them as polynomials in $z, \ \overline z$.
  For this reason  the new notation $Z^{n,k}$ is introduced
  according to the transformation of indexes $l=n-2k$ in (\ref{zernike}).
  Thus we get
 \begin{equation}
 V_{n,n-2k}(r,\psi )=
 z^{n-2k} P_k^{(0,|n-2k|)}(2|z|^2-1), \ (k=0,1,...,n). \nonumber
 \nonumber
 \end{equation}
 So, we define
 \begin{equation}
 \label{zecomp}
 Z^{n,k}(z,\overline{z}):=(-1)^k V_{n,n-2k}(r,\psi)=(-1)^kz^{n-2k} P_k^{(0,|n-2k|)}(2|z|^2-1),
 \end{equation}
 where $n=0,1,2,...$ and $k=0,1,...,n.$
 The first index (superscript) $n$ indicates the degree of a polynomial $Z^{n,k}$ and the
 second superscript $k$ denotes its order in a bunch
 $Z^{n,0},Z^{n,1},...,Z^{n,n}.$
 The multiplier $(-1)^k$ in (\ref{zecomp})
 was introduced for convenience of further computations.

 This definition can also be rewritten as
 \begin{equation}
 \label{zer1}
 Z^{n,k}(z,\overline{z})=
  \begin{cases}
 \sum _{s=0}^{k} C^s_kC^s_{n-k}z^{n-k-s}
 (1- z \overline z)^s (- \overline z)^{k-s}   &\text{for}
 \ k=0,1,..., \left [  \dfrac{n}{2} \right ]
 \\
 \\
 (-1)^n {\overline Z}^{n,n-k}(z,\overline{z}) &\text{for} \ k=\left [  \dfrac{n}{2}
 \right ]+1,..., n,\\
 \end{cases}
 \end{equation}
 where $[\cdot]$ denotes the integer part of a number.
 \par
 After the evaluation of (\ref{zer1}) we get
 \begin{equation}
 \label{zer2}
 Z^{n,k}(z,\overline{z})=
  \begin{cases}
 \sum _{s=0}^{k} (-1)^{k-s}C^s_{n-k}C^{k-s}_{n-s} \overline z ^{k-s}
 z^{n-k-s} &\text{for} \ k=0,1,..., \left [  \dfrac{n}{2} \right ]\\
  \\
 (-1)^n {\overline Z}^{n,n-k}(z,\overline{z}) &\text{for} \ k=\left [  \dfrac{n}{2} \right ]+1,..., n.\\
 \end{cases}
 \end{equation}
 For example, the first Zernike polynomials up to the degree (order) $n=4$ are
 \par
 {\footnotesize
 $$\begin{array}{lllll}
 Z^{0,0}=1 & & & & \\
 \\
 Z^{1,0}=z &  Z^{1,1}=-\overline z & & &\\
 \\
 Z^{2,0}=z^2 & Z^{2,1}=1-2z\overline z & Z^{2,2}=\overline z^2 & & \\
 \\
 Z^{3,0}=z^3 & Z^{3,1}=2z-3z^2\overline z & Z^{3,2}=3z\overline z^2-2\overline z &
 Z^{3,3}=-\overline z^3 & \\
 \\
 Z^{4,0}=z^4 & Z^{4,1}=3z^2-4z^3\overline z & Z^{4,2}=1-6z\overline z+6z^2\overline z^2 &
  Z^{4,3}=3\overline z^2-4z\overline z^3 &
 Z^{4,4}=\overline z^4.
 \end{array}
 $$
 }
 \\
 \par The Zernike polynomials are orthogonal in the unit disc $\mathbb D$,
 obey the following orthogonality relation
 \begin{equation}
 \label{ortzer}
 \langle\langle Z^{n,k}, ~Z^{m,s}\rangle\rangle=\frac {\pi}{n+1}\delta_{n,m}\delta_{k,s}
 \end{equation}
 and their $L_2$-norms are equal to
 \begin{equation}
 \nonumber
 ||Z^{n,k}||=\sqrt{\frac {\pi}{n+1}}, \  \  (k=0,1,...,n) .
 \end{equation}
 \par
 It allows the expansion  of an arbitrary function $a (z,\overline z) \in L_2({\mathbb D})$
 in terms of a unique combination of Zernike polynomials.
 \begin{equation}
 \label{zerdec}
 a(z,\overline z)=\sum_{n=0}^{\infty}\frac {n+1}{\pi}
 \sum_{k=0}^{n}\langle\langle a,~Z^{n,k}\rangle\rangle Z^{n,k}(z,\overline z).
 \end{equation}
 \par  Since  we use the complex variables $z$ and $\overline z$ and treat
 them as independent variables here, we'll sometimes write $a(z)$
 instead of $a(z,\overline z).$ The formal partial derivatives with respect to
 $z$ and $\overline z$ are defined in the usual way by (\ref {part}).
 \par  In the following theorem we formulate in complex variables some novel properties
 of Zernike polynomials.
 \par {\theorem 1 The following properties take place:
 \par \noindent (a) Zernike polynomials (\ref{zer1})
  have the differential representation
 \begin{equation}
 \label{a}
 Z^{n,k}(z,\overline{z})= \frac {1}{k!} \frac {\partial ^k}
 {\partial z^{k}}
 \left [ z^{n}\left (\frac{1}{z}-\overline{z} \right )^{k}
 \right ], \ (n\geq 0, \  k=0,1,...,n).
 \end{equation}
 \par \noindent (b) Zernike polynomials (\ref{zer1})
 are the solution
 of the elliptic system
 \begin{eqnarray}
 \label{c}
 \left \{
 \begin{array} {lll}
                         &(Z^{n,n})_z&=\ 0 \\
 (Z^{n,n})_{\overline{z}}\ \ +&(Z^{n,n-1})_z&=\ 0  \\
 &...&  \\
 (Z^{n,k})_{\overline{z}}\ \ +&(Z^{n,k-1})_z&=\ 0 \\
 &...&  \\
 (Z^{n,1})_{\overline{z}}\ \ +&(Z^{n,0})_z&=\ 0 \\
 (Z^{n,0})_{\overline{z}}& &=\ 0\\
 \end{array} \right.
 \end{eqnarray}
 and satisfy  boundary conditions
 \begin{eqnarray}
 \label{cc}
 Z^{n,k}(t,\overline{t})= (-1)^kt^{n-2k}, \ |t|=1, \ (n\geq 0, \
 k=0,1,...,n).
 \end{eqnarray}
 \par \noindent (c) Zernike polynomials (\ref{zer1}) can be represented in the form of Cauchy-type
 integral
 \begin{equation}
 \label{b}
 \frac {1}{2\pi {\i}} \mathop{\mathop {\int}}\limits_{|t|=1}
 \frac {t^{n}(\overline t-\overline z)^k} {(t-z)^{k+1}}{\d}t=
  \begin{cases}
 Z^{n,k}(z,\overline{z}) \ &\text{for} \ n \geq 0, \ k=0,1,...,n \\
 \\
 0                         &\text{for} \ n \geq 0, \ k>n \  \text{or}  \ k<0.
 \end{cases}
 \end{equation}
 }
 {\proof
 (a) Let's first prove (\ref{a}) for $ k=0,1,..., \left [ \dfrac{n}{2} \right ]$.
 \par \noindent By Leibnitz formula $\Big [ uv \Big ]^{(k)}=\sum_{s=0}^{k}C_k^{s}u^{(s)}v^{(k-s)}$
 we get
 \begin{eqnarray}
 \left [ z^n\left (\frac {1}{z}-\overline z\right )^k \right ]^{(k)}_{z}=
 \Big  [  z^{n-k}(1- z \overline z)^k \Big  ]^{(k)}_{z}= \sum _{s=0}^{k}
 C^s_k \Big [z^{n-k} \Big ]^{(s)}_{z}\Big [(1- z \overline z)^k \Big ]^{(k-s)}_z \nonumber
 \end{eqnarray}
 \begin{eqnarray} \label{a1}
 =&&\sum _{s=0}^{k}
 C^s_k(n-k)(n-k-1)...(n-k-s+1)z^{n-k-s}\nonumber \\
 && \times  k(k-1)...(s+1) (1- z \overline z)^s (-\overline z)^{k-s}\nonumber \\
 =&&k!\sum _{s=0}^{k} C^s_kC^s_{n-k}z^{n-k-s}
 (1- z \overline z)^s (-\overline z)^{k-s}=  k!Z^{n,k}(z,\overline{z}).
 \end{eqnarray}
 \par \noindent Now let's substitute $k\to n-k$ in (\ref{a1}). Taking into account that $k\leq
 n-k,$ we get
 \begin{eqnarray}
 \label{a2}
 &&\left [ z^n \left (\frac {1}{z}- \overline z \right )^{n-k} \right ]^{(n-k)}_{ z}
 = \Big  [ z^{k}(1- z \overline z)^{n-k} \Big  ]^{(n-k)}_{ z}\nonumber\\
 =&& \sum_{s=0}^{k} C^s_{n-k}\Big [z^{k} \Big ]^{(s)}_{ z} \Big [(1- z \overline z)^{n-k}\Big ]^{(n-k-s)}_{z}
 \nonumber \\
 =&&\sum _{s=0}^{k} C^s_{n-k}k(k-1)...(k-s+1) z^{k-s} \nonumber \\
  && \times (n-k)(n-k-1)...(s+1)(1- z \overline z)^s (- \overline z)^{n-k-s}  \nonumber \\
  =&&(n-k)!\sum _{s=0}^{k} C^s_kC^s_{n-k}(-1)^n(\overline z)^{n-k-s}
 (1- z \overline z)^s (- z)^{k-s}\nonumber\\
 =&&(-1)^n (n-k)!\overline Z^{n,k}(z,\overline{z})= (n-k)!Z^{n,n-k}(z,\overline{z}).
 \nonumber
 \end{eqnarray}
 And the assertion (a) follows.
 \par \noindent (b) It is easy to verify equations (\ref{c}) by direct computation
  using formula
 (\ref{zer2}). On the boundary of the unit disc $\mathbb D,$
 due to the normalization  $P_k^{(0,|n-2k|)}(1)=1$ of  Jacobi polynomials,
 we get
 \begin{eqnarray}
 \label{bouzer}
 Z^{n,k}(t,\overline{t})= (-1)^kt^{n-2k}, \ |t|=1, \ (n\geq 0, \  k=0,1,...,n).
 \nonumber
 \end{eqnarray}
 \par \noindent(c) In (\ref{b}) we take advantage of Newtonian binomial formula
 \begin{eqnarray}
 \label{b1}
 &&\frac {1}{2\pi {\i}}
 \mathop {\int}\limits_{|t|=1} \frac {t^{n}(\overline t-\overline z)^k}
 {(t-z)^{k+1}}{\d}t=\sum_{s=0}^{k}C_k^s(-\overline{z})^{k-s}
 \frac {1}{2\pi {\i}}
 \mathop {\int}\limits_{|t|=1} \frac {t^{n}{\overline t}^s}
 {(t-z)^{k+1}}{\d}t  \nonumber \\
 &&=\sum_{s=0}^{k}C_k^s(-\overline{z})^{k-s}
 \frac {1}{2\pi {\i}}
 \mathop {\int}\limits_{|t|=1} \frac {t^{n-s}}
 {(t-z)^{k+1}}{\d}t \nonumber \\
 &&=\sum_{s=0}^{k}C_k^s(-\overline{z})^{k-s}\frac {1}{k!}
 \frac {k!}{2\pi {\i}}
 \mathop {\int}\limits_{|t|=1} \frac {t^{n-s}}
 {(t-z)^{k+1}}{\d}t  \nonumber \\
 &&=  \sum_{s=0}^{k}C_k^s(-\overline{z})^{k-s}\frac {1}{k!}
 \frac {{\d}^k}{{\d}z^k}\Big [ z^{n-s} \Big ]=
 \frac {1}{k!}
 \frac {\partial ^k}{\partial z^k}\left [
 \sum_{s=0}^{k}C_k^s(-\overline{z})^{k-s}
 z^{n-s} \right ]  \nonumber \\
 &&= \frac {1}{k!}
 \frac {\partial ^k}{ \partial z^k}\left [
 z^{n}\left (\frac{1}{z}-\overline{z} \right )^{k}
 \right ]. \nonumber
 \end{eqnarray}
 \noindent Theorem 1 is proved.
 \endproof

 \subsection {Fan-beam Radon transform of Zernike polynomials}
 The fan-beam Radon transform ${\cal D}$ of a scalar
 function $a(x^1,x^2)$ is defined by (\ref{comfanT}) for $m=0$.
 We have
 \begin{equation}\label{scfanT}
 [{\cal D}{a}](\beta,\varphi)
 =\mathop{\int}\limits_{\tau(t,\varphi)}^{\quad t}{a}(\zeta,\overline \zeta)|{\d}\zeta|,
 \end{equation}
 where $ t=e^{{\i}\beta}, \ \beta \in [0,2\pi), \ \tau(t,\varphi)=-\overline te^{2{\i}\varphi},
 \ |\beta-\varphi |\leq \dfrac {\pi}{2}, \ \varphi= \arg(t-\tau).$
 For  $|\beta-\varphi |> \frac {\pi}{2}$ we complete the definition of the
 fan-beam transform
 (\ref {scfanT}) with the condition
 $$[{\cal D}{a}](\beta,\varphi)=-[{\cal D}{a}](\beta,\varphi+\pi).$$

 {\theorem 2 The fan-beam Radon transform ${\cal D}Z^{n,k}$
 of Zernike polynomials $Z^{n,k}$ equals to
 \begin{equation}
 \label{e}
 [{\cal D}Z^{n,k}](\beta,\varphi)=
 \frac {2e^{{\i}(n-2k)\varphi}}{n+1} \times
 \begin{cases}
 \cos[(n+1)(\beta-\varphi)] \  \ \ \text{for} \ n=\text{even} \\
 \\
 {\i}\sin[(n+1)(\beta-\varphi)] \  \ \text{for} \ n=\text{odd,}\\
 \end{cases}
 \end{equation}
 where $\beta \in [0,2\pi), \ \varphi \in [0,2\pi).$
 }

 \noindent  \proof
 At first we  introduce  the auxiliary polynomials $X^{n,k}$ defined by
 $$
 X^{n,k}(z,\overline z):= \frac {1}{k!} \frac {\partial ^{k-1}}
   {\partial z^{k-1}}
   \left [ z^{n}\left (\frac{1}{z}-\overline{z}
   \right )^{k}
   \right ], \    (n \geq 1, \ k=1,...,n).
 $$
 Then the next equations  follows directly  from (\ref{a})
 \begin{equation}\label{partial}
 \frac {\partial X^{n,k}}{\partial z}
 =Z^{n,k}(z,\overline{z}), \ \frac {\partial X^{n,k}}{\partial \overline z}
 =-Z^{n,k-1}(z,\overline{z}).
  \end{equation}
 Using  (\ref{a}) we can  verify, that
 \begin{equation}
 X^{n,k}=\frac {1}{k} \left ( Z^{n-1,k-1}-\overline{z}Z^{n,k-1} \right ).
 \nonumber
 \end{equation}
 Then combining above and  (\ref{cc}) we obtain the  boundary conditions
 \begin{equation}
 \label{xgran}
 X^{n,k}(t,\overline{t}) =0, \ |t|=1.
 \end{equation}
  Let's compute the fan-beam transformation of Zernike polynomials.
 To this end  we use previously obtained derivatives (\ref{partial})
  for computing by (\ref{derivteta}) the derivative of $X^{n,k}$ in the direction
 ${\bf \Theta}=\begin{pmatrix} \Theta^1  \\ \Theta^2 \end{pmatrix}
 = \begin{pmatrix}e^{{\i}\varphi} \\e^{-{\i}\varphi} \end{pmatrix}.$
 So we get
 \begin{eqnarray}
 \frac {\partial X^{n,k}}{\partial{\bf \Theta}}=
 e^{{\i}\varphi}\frac {\partial X^{n,k}}{\partial z}+
 e^{-{\i}\varphi}\frac {\partial X^{n,k}}{\partial \overline z}=
 e^{{\i}\varphi}Z^{n,k}-e^{-{\i}\varphi}Z^{n,k-1},
 \nonumber
 \end{eqnarray}
 or the same in another form
 \begin{eqnarray}
  Z^{n,k}(z,\overline{z})=e^{-2{\i}\varphi}Z^{n,k-1}+e^{-{\i}\varphi}
  \frac {\partial X^{n,k}}{\partial{\bf \Theta}}. \nonumber
\end{eqnarray}
 This equation combined with (\ref{xgran}) is used for evaluation of the next integral
 \begin{eqnarray}
 \mathop {\int}\limits_{\tau}^{\hspace{2mm} t}
 Z^{n,k}|{\d}\zeta|&=& e^{-2{\i}\varphi}
 \mathop {\int}\limits_{\tau}^{\hspace{2mm}t}
 Z^{n,k-1}|{\d}\zeta| +e^{-{\i}\varphi} (X^{n,k}(t,\overline t)-X^{n,k}(\tau,\overline \tau) )
 \nonumber \\
 &=& e^{-2{\i}\varphi}
 \mathop {\int}\limits_{\tau}^{\hspace{2mm}t} Z^{n,k-1}|{\d}\zeta|,
  \nonumber
 \end{eqnarray}
 where  $t=e^{{\i}\beta}, \ \tau(t,\varphi)=-\overline te^{2{\i}\varphi}.$
 Unwrapping the recurrence relation gives
 $$\mathop{\mathop {\int}}\limits_{\tau}^{\quad t}Z^{n,k}|{\d}\zeta|=
 e^{-2k{\i}\varphi} \mathop{\mathop {\int}}\limits_{\tau}^{\quad t}
 Z^{n,0}|{\d}\zeta|.
$$
 The last integral is computed directly, taking into account that
 $Z^{n,0}(z,\overline z )=z^n$
 \begin{eqnarray}
 &&\mathop{\int}\limits_{ \tau}^{\quad t}
 Z^{n,0}|{\d}\zeta|=
 \mathop {\int}\limits_{ \tau}^{\quad t}\zeta^n |{\d}\zeta|=
  \mathop {\int}_{0}^{|t-\tau|}( \tau +se^{{\i}\varphi})^n{\d}s \nonumber \\
  &=&  \frac {1}{n+1} ( \tau +se^{{\i}\varphi})^{n+1}e^{-{\i}\varphi}
  \Big |^{|t-\tau|}_0 = \frac
 {1}{n+1}e^{-{\i}\varphi }(t^{n+1} - \tau^{n+1}). \nonumber
  \end{eqnarray}
 Finally, we get
 \begin{eqnarray}[{\cal D}Z^{n,k}](\beta,\varphi)&=&
 \mathop{\int}\limits_{\tau}^{\hspace{2mm} t}
 Z^{n,k}|{\d}\zeta|=
 \frac {e^{-2k{\i}\varphi}}{n+1}(e^{-{\i}\varphi}t^{n+1}+(-1)^{n} e^{{\i}(2n+1)\varphi}\overline{t}^{n+1})
 \nonumber \\
 &=& \frac{2e^{{\i} (n-2k)\varphi}}{n+1} \times
 \begin{cases}
 \cos[(n+1)(\beta-\varphi)] \  \ \ \text{for} \ n=\text{even} \\
 \\
 {\i}\sin[(n+1)(\beta-\varphi)] \  \ \text{for} \ n=\text{odd}.
 \end{cases}
 \nonumber
 \end{eqnarray}
 Theorem 2 is proof.
 \endproof

\section { Construction of the orthogonal polynomial basis
and SVD}

 In this section we describe the construction of  orthogonal polynomial basis in
 the space of solenoidal (divergence free) tensor fields ${\bf H}(\mathbb D;{\bf S}_m,\delta=0).$

 We are given a polynomial of degree  $N$ solenoidal  $m$-covariant
 tensor field ${\bf a} \in {\bf H}_N(\mathbb D;{\bf S}_m,\delta=0)$
 and in complex coordinates we have ${\bf a} \rightarrowtail {\bf A}=\{A_k \}.$
 As was mentioned earlier, we use a pseudovector notation for the tensor field ${\bf A}$
 \begin{equation}
 {\bf A}=\begin{pmatrix} A_m  \\ A_{m-1}\\ ...\\ A_1\\ A_0
 \end{pmatrix}.
 \label{A}
 \end{equation}
 The condition (\ref{cond}) now looks like
 \begin{equation}
 \begin{pmatrix} A_m  \\ A_{m-1}\\ ...\\ A_1\\ A_0
 \end{pmatrix}=
 \begin{pmatrix} \overline A_0  \\
 \overline A_{1}\\ ...\\ \overline A_{m-1}\\ \overline A_{m}
 \end{pmatrix}.
 \label{A1}
 \end{equation}
 For  $m,\ n \geq 0, \ k=0,...,n+m$ and  $2k \neq m+n$ we define
 polynomial of degree $n$ symmetric tensor fields (in complex variables )
 \begin{equation}
 {\bf S}^{(+m)}_{n,k}:=(-1)^n
 \begin{pmatrix}
 Z^{n,k} +    \overline Z^{n,k-m} \\ Z^{n,k-1} +  \overline Z^{n,k-m+1}\\ ...\\ Z^{n,k-m+1}
 +\overline Z^{n,k-1}\\ Z^{n,k-m}+ \overline Z^{n,k}
 \end{pmatrix},
 \ {\bf S}^{(-m)}_{n,k}:=\dfrac{1}{{\i}} \begin{pmatrix}
 Z^{n,k} -    \overline Z^{n,k-m} \\ Z^{n,k-1} -  \overline Z^{n,k-m+1}\\ ...\\
 Z^{n,k-m+1}-\overline Z^{n,k-1}\\ Z^{n,k-m}- \overline Z^{n,k}
 \end{pmatrix},
 \label{Sbas}
 \end{equation}
 where for convenience we set $Z^{n,k} \equiv 0$ for $k<0$ or $k>n.$
 \par For $2k = m+n$ we have two  cases:
 \\
 \noindent the first one is when $m$ and $n=$  even
 \begin{equation}
 {\bf S}^{(+m)}_{n,\frac {m+n}{2}}:=
 \begin{pmatrix} Z^{n,\frac {m+n}{2}} \\ Z^{n,\frac {m+n}{2}-1} \\ ...
 \\ Z^{n,\frac {m+n}{2}-m+1} \\ Z^{n,\frac {m+n}{2}-m}
 \end{pmatrix}, \
 {\bf S}^{(-m)}_{n,\frac {m+n}{2}}:=0;
 \label{S1}
 \end{equation}
 the second one is when $m$ and $n=$  odd
   \begin{equation}   {\bf S}^{(+ m)}_{n,\frac {m+n}{2}}:=0, \
 {\bf S}^{(-m)}_{n,\frac {m+n}{2}}:=
 \dfrac{1}{{\i}} \begin{pmatrix} Z^{n,\frac {m+n}{2}} \\ Z^{n,\frac {m+n}{2}-1}
 \\ ...\\ Z^{n,\frac {m+n}{2}-m+1}\\ Z^{n,\frac {m+n}{2}-m}
 \end{pmatrix}.
 \label{S2}
 \end{equation}
 \\
 \par  {\bf Remark 1.} Polynomial tensor fields (\ref{S1}) and  (\ref{S2})
  can be evaluated by general formulae (\ref{Sbas}), but then
  the result should be divided by 2.
  \\
 \par  {\bf Remark 2.} Note, that for  $k+s=m+n$ the equation
 $${\bf S}^{(\pm m)}_{n,k}=(-1)^n{\bf S}^{(\pm m)}_{n,s}
 $$
 takes place, therefore in (\ref{Sbas}) we set only
 $k=0,1,...,\left [\frac {n+m}{2} \right ],$
  where $[\cdot]$  defines the integer part of a number.
 \par

 So, if we now make transformations (\ref{rtocom}) or (\ref{a_k})
 from complex variables   to real variables
 \begin{equation}\label{sbas}
 {\bf S}^{(\pm m)}_{n,k}(z, \overline z) \to {\bf s}^{(\pm m)}_{n,k}(x^1, x^2), \
 \left (n \geq 0, \ k=0,1,...,\left [\frac {n+m}{2} \right ]\right ),
 \end{equation}
 then we get  polynomial real-valued tensors ${\bf s}^{(\pm m)}_{n,k}$
 in Cartesian variables $(x^1, x^2).$
 \par
 {\lemma 1 The tensor fields
 ${\bf s}^{(\pm m)}_{n,k} \ \left (n=0,...,N, \ k=0,1,..., \left [\frac{n+m}{2} \right ]
  \right )$ \\
  defined by
 (\ref{sbas})
 form an orthogonal basis of finite-dimensional  subspace
 \\ ${\bf H}_N(\mathbb D;{\bf S}_m,\delta=0),$  thus
 $$\dim {\bf H}_N(\mathbb D;{\bf S}_m,\delta=0)
  =\dfrac {(N+1)(N+2+2m)}{2}.$$ }
 \par \proof
 Consider a tensor field ${\bf a}\in {\bf H}_N(\mathbb D;{\bf S}_m,\delta=0)$
 and let ${\bf a} \rightarrowtail {\bf A}.$
 Expand each component  $A_k$
 of pseudovector (\ref{A}) in the sum
 of  Zernike polynomials and  use the solenoidality condition (\ref{sol}).
 Taking into account the property (b) from Theorem 1,
 we get an expansion of tensor ${\bf A}$ into the sum
 \begin{equation}
 {\bf A}=\begin{pmatrix} A_m  \\
 A_{m-1}\\ ...\\ A_1\\ A_0
 \end{pmatrix}
 =\sum_{ n=0}^{N}\sum_{k=0}^{n+m}c_{n,k}
 \begin{pmatrix}
 Z^{n,k}  \\ Z^{n,k-1}\\ ...\\ Z^{n,k-m+1}\\ Z^{n,k-m}
 \end{pmatrix}.
 \label{sum}
 \end{equation}
 From another hand, (\ref{A1}) yields
 \begin{equation}
 {\bf A}= \begin{pmatrix}
 \overline A_{0}\\
 \overline A_{1}\\
 ...\\
 \overline A_{m-1}\\
 \overline A_{m}
 \end{pmatrix}=
 \sum_{n=0}^{N}\sum_{k=0}^{n+m}\overline c_{n,k}
 \begin{pmatrix}
 \overline Z^{n,k-m}  \\ \overline Z^{n,k-m+1}\\ ...\\ \overline Z^{n,k-1}\\ \overline Z^{n,k}
 \end{pmatrix}=
 \sum_{n=0}^N \sum_{k=0}^{n+m}\overline c_{n,n+m-k}
 \begin{pmatrix}
 \overline Z^{n,n-k}  \\ \overline Z^{n,n-k+1}\\ ...\\ \overline Z^{n,n-k+m-1}\\ \overline
 Z^{n,n-k+m}
 \end{pmatrix} \nonumber
 \end{equation}
 \begin{equation}
  =\sum_{n=0}^N\sum_{k=0}^{n+m}(-1)^n\overline c_{n,n+m-k}
 \begin{pmatrix}
 Z^{n,k}  \\ Z^{n,k-1}\\ ...\\ Z^{n,k-m+1}\\ Z^{n,k-m}
 \end{pmatrix}.
 \end{equation}
 Comparing the last expression with (\ref{sum}), we obtain
 \begin{equation}c_{n,k}=(-1)^n\overline c_{n,n+m-k}, \
 (k=0,...,n+m) . \label {coef}
\end{equation}
 Splitting the coefficients $c_{n,k}$ in (\ref{sum}) into the real and imaginary parts
  $$c_{n,k}=a_{n,k}+{\i}b_{n,k}, \   \left (k=0,...,\left [\frac {n+m}{2} \right ] \right ),
  $$
 taking into account (\ref{coef}) and definition of ${\bf S}^{\pm m)}_{n,k},$ we  finally
 get
  $$\begin{pmatrix} A_m  \\ A_{m-1}\\ ...\\ A_1\\ A_0
 \end{pmatrix}=
 \sum_{n=0}^{N} \sum_{k=0}^{[\frac {n+m}{2}]*}a_{n,k}
 \begin{pmatrix}
 Z^{n,k} +    \overline Z^{n,k-m} \\ Z^{n,k-1} +
 \overline Z^{n,k-m+1}\\ ...\\ Z^{n,k-m+1}
 +\overline Z^{n,k-1}\\ Z^{n,k-m}+   \overline Z^{n,k}
 \end{pmatrix}
 + {\i} b_{n,k}\begin{pmatrix} Z^{n,k} -    \overline Z^{n,k-m} \\
 Z^{n,k-1} -  \overline Z^{n,k-m+1}\\
 ...\\ Z^{n,k-m+1}-\overline Z^{n,k-1}\\ Z^{n,k-m}-   \overline Z^{n,k}
 \end{pmatrix} \nonumber
 $$
 \begin{eqnarray} =\sum_{n=0}^N\sum_{k=0}^{[\frac {n+m}{2}] \ *}
 (-1)^na_{n,k}{\bf S}^{(+m)}_{n,k} (z,\overline z)
 -b_{n,k}{\bf S}^{(-m)}_{n,k} (z,\overline z). \nonumber \label{4.8}
 \end{eqnarray}
 The sign $*$ here means that in the case of even $n$ and even $m$
 the coefficient $b_{n,\frac{n+m}{2}}$ should be set to  $0$,
 and in the case of odd $n$ and odd $m$ the coefficient $a_{n,\frac {n+m}{2}}$ should be set to $0.$
 \par
 So, we have that  the tensor field ${\bf a}$ is a linear combination of polynomial
 tensor fields (\ref {sbas}).
 \par
 Now we show that polynomial tensor fields
 (\ref {sbas}) are orthogonal. Let $k \neq s, \ k,s=0,1,...,\left [\frac {n+m}{2} \right ]$
 and remark, that  then  $k+s \neq m+n$ take place.
 Using formula (\ref {pscal}) we have
 \begin{eqnarray}
  &&\langle\langle {\bf s}^{(\pm m)}_{n,k},~{\bf s}^{(\pm m)}_{n,s}\rangle  \rangle=
  \langle\langle {\bf S}^{(\pm m)}_{n,k},~{\bf S}^{(\pm m)}_{n,s}\rangle \rangle  \nonumber \\
 &=& \pm  2^m \iint\limits_{\mathbb D}\sum_{p=0}^{m}C^p_m(Z^{n,k-m+p}\pm \overline Z^{n,k-p})
 (Z^{n,s-p}\pm \overline Z^{n,s-m+p}){\d}V^2.  \nonumber
 \end{eqnarray}
  Taking into account the orthogonality  of Zernike polynomials
  we obtain
 \begin{eqnarray}
  \langle\langle {\bf s}^{(\pm m)}_{n,k},~{\bf s}^{(\pm m)}_{n,s}\rangle \rangle =0,
  \ k\neq s. \nonumber
  \end{eqnarray}
 Let's now evaluate the norms of polynomial tensor ${\bf s}^{(\pm m)}_{n,k}$
  by the formula (\ref {norm}).
 At first we consider the case $n+m \neq 2k,$ then
 \begin{eqnarray}
  ||{\bf s}^{(\pm m)}_{n,k}||^2&=&
 ||{\bf S}^{(\pm m)}_{n,k}||^2= 2^m \iint\limits_{\mathbb D}
 \sum_{p=0}^{m}C^p_m|Z^{n,k-m+p}\pm \overline Z^{n,k-p}|^2{\d}V^2 \nonumber \\
 &=&2^m \sum_{p=0}^{m}C^p_m
 \left (||Z^{n,k-m+p}\pm \overline Z^{n,k-p}||^2 \right )\nonumber \\
 &=&2^m \sum_{p=0}^{m}C^p_m||Z^{n,k-m+p}||^2+|| \overline Z^{n,k-p}||^2 \nonumber \\
 &=& 2^{m+1}\sum_{p=k-n}^{k}C^p_m||Z^{n,k-p}||^2 =
 \frac {2^{m+1} \pi}{n+1} \alpha^{(m)}_{n,k},
 \label{alph}
 \end{eqnarray}
 where coefficients $\alpha^{(m)}_{n,k}$ are defined by
 \begin{equation}\label{alph1}
 \alpha^{(m)}_{n,k}=\sum^{k}_{p=k-n}C^p_m, \ \left ( n\geq 0, \
 k=0,1,..., \left [\dfrac {m+n}{2} \right] \right ).
 \end{equation}
\par {\remark 3 In this formula for convenience we set $C^p_m=0$ if $p<0$ or  $p>m.$}
\\
 If  $n+m = 2k,$  i.e.   $k=\frac {n+m}{2},$
 then taking into account the above calculations, we get
 $$ ||{\bf s}^{(+ m)}_{n,\frac {m+n}{2}}||^2=
 \begin{cases}
 0   &\text {for \ n=odd} \\
 \\
 \dfrac {2^{m}\pi}{n+1}\alpha^{(m)}_{n,\frac {m+n}{2}}$$
 \ &\text{for \ n=even},
 \end{cases}
 $$
 $$
 ||{\bf s}^{(-m)}_{n,\frac {m+n}{2}}||^2=
 \begin{cases}
 \dfrac {2^{m}\pi}{n+1}\alpha^{(m)}_{n,\frac {m+n}{2}}  \ &\text {for \ n=odd} \\
 \\
 0 \ &\text {for \ n=even.}
 \end{cases}
  $$
 Obviously,
 $$\dim{\bf H }_{N}(\mathbb D;{\bf S}_m,\delta=0)  =\sum_{n=m}^{m+N}n=\frac
{(N+1)(N+2+2m)}{2}.$$
 Lemma 1 is proved.
 \endproof
 \par Lemma 1 and
 the definition of the  subspace of solenoidal tensor fields
 ${\bf H}(\mathbb D;{\bf S}_m, \delta=0)$ yield
 \par
 {\corollary 1
  Polynomial tensor fields
 ${\bf s}^{(\pm m)}_{n,k} \ \left (n\geq 0, \ k=0,1,...,\left
 [\frac{n+m}{2} \right ] \right )$ \\
 form an orthogonal basis in the subspace of solenoidal tensor fields \\
 ${\bf H}(\mathbb D;{\bf S}_m, \delta=0)   \subset {\mathbf L}_2( \mathbb D ; {\bf S}_m).$
 }

  Now we are in a position to define the  SVD for the  fan-beam Radon transform
  ${\cal D}_m.$

 {\theorem 3 The singular values of the operator $(\ref{fanT})$
 $${\cal D}_m:{\mathbf L}_{2}(\mathbb D; {\bf S}_m ) \to L_2([0,2\pi)\times[0,2\pi))
 $$
 are given by
 \begin{equation}\label{sval} \sigma^{(m)}_{n,k}\equiv
 \sigma^{(\pm m)}_{n,k}:=
 \sqrt {\frac {8\pi }{(n+1)2^{m}} \alpha^{(m)}_{n,k}},
 \left  (n\geq 0, \ k=0, 1,..., \left [\frac {n+m}{2} \right ] \right),
 \nonumber
 \end{equation}
 where coefficients  $\alpha^{(m)}_{n,k}$
 are defined by the formula (\ref {alph1}).
 If a solenoidal real-valued symmetrical tensor field
 ${\bf a}(x^1,x^2) \in {\mathbf  L}_2(\mathbb D; {\bf S}_m)$
 has an expansion
 \begin{equation}\label{a=}
 {\bf a}(x^1,x^2)=
 \sum_{n=0}^{\infty } \sum_{k=0}^{[\frac {n+m}{2}] \ \ast}\frac {1}{
 ||{\bf s}^{(\pm m)}_{n,k}|| } \left (a_{n,k}{\bf s}^{(+m)}_{n,k} (x^1,x^2)+
 b_{n,k}{\bf s}^{(-m)}_{n,k}
 (x^1,x^2) \right ),
 \end{equation}
 then the fan-beam Radon transform ${\cal D}_m{\bf a}$ has the following singular value
 decomposition
 \begin{equation}
 [{\cal D}_m{\bf a}](\beta,\varphi)=\sum_{n=0}^{\infty}
 \sum_{k=0}^{[\frac {n+m}{2}]
 \ \ast}\sigma^{(m)}_{n,k}\left (a_{n,k}f^{(+m)}_{n,k}(\beta,\varphi)
 +b_{n,k}f^{(-m)}_{n,k}(\beta,\varphi) \right ),
 \label{Da}
 \end{equation}
 where singular functions are
 $$f^{(+m)}_{n,k}(\beta,\varphi):=\frac {1}{\pi}
 \begin{cases}
 \cos[(n+1)(\beta-\varphi)]\cos[(n-2k+m)\varphi] \\
 \\
 \sin[(n+1)(\beta-\varphi)]\sin[(n-2k+m)\varphi], \\
 \end{cases}
 $$
 $$ f^{(-m)}_{n,k}(\beta,\varphi): = \frac {1}{\pi}
 \begin{cases}
 \cos[(n+1)(\beta-\varphi)]\sin[(n-2k+m)\varphi] \\
 \\
 \sin[(n+1)(\beta-\varphi)]\cos[(n-2k+m)\varphi], \\
 \end{cases}
 $$
 when $n \geq 0, \ k=0,..., \left [ \frac {n+m}{2} \right ]$
 and $2k \neq m+n$ and
 $$ f^{(+m)}_{n,\frac {m+n}{2}}(\beta,\varphi):=\frac {1}{\sqrt 2\pi}
 \begin{cases}
 \cos[(n+1)(\beta-\varphi)] \\
  0,
 \end{cases}
 $$
 $$ f^{(-m)}_{n,\frac {m+n}{2}}(\beta,\varphi): = \frac {1}{\sqrt 2 \pi}
 \begin{cases}
 0 \\
 \sin[(n+1)(\beta-\varphi)], \\

 \end{cases}
 $$
 when $2k=m+n.$
 In all expressions above top line corresponds to the even values of $n$,
 and bottom line --- to the odd $n$.
 The sign $*$  in (\ref{a=}) and (\ref{Da}) near by the inner sum denotes
 that in the case of even $n$ and $m$
 the coefficient $b_{n,\frac {n+m}{2}}$ should be set to $0,$
 and in the case of odd $n$ and $m$ ---
 the coefficient $a_{n,\frac {n+m}{2}}$ respectively.
 }
 \par \proof Note, that the system of functions
 $f^{(\pm  )}_{n,k}$ for $n\geq 0,  \ k=0,1,...,\left [\frac {n+m}{2}\right ]$
 is the subsystem of the standard orthonormal
 basis of  $L_2([0,2\pi)\times[0,2\pi))$  and there is
 the basis of the image of the tensor fan beam transform  ${\cal D}_m$.
 \par Let's evaluate now the fan-beam transform ${\cal D}_m$ for basis polynomial
 tensor (\ref{sbas}). For this we introduce
 $$
 {\bf A}:=\begin{pmatrix} Z^{n,k}  \\  Z^{n,k-1}\\ ...\\  Z^{n,k-m+1}\\  Z^{n,k-m}
 \end{pmatrix}, \ A_{m-s}=Z^{n,k-s}
 \ \ \text{and} \ \
 {\bf B}:=\begin{pmatrix}
 \overline Z^{n,k-m}
 \\ \overline Z^{n,k-m+1}\\ ...\\ \overline Z^{n,k-1}\\ \overline
 Z^{n,k}
 \end{pmatrix}, \ B_s=\overline Z^{n,k-s}.
 $$
 Hence in case of $n+m\neq 2k$ we get
 $$
 {\bf S}^{(+m)}_{n,k}=(-1)^n({\bf A}+{\bf B}), \
 {\bf S}^{(-m)}_{n,k}=\frac{1}{{\i}}({\bf A}-{\bf B}).
 $$
 \noindent  Then  by using formulae (\ref{Vtr1}) and Theorem 2 we have
 \begin{eqnarray}
 && [{\cal D}_m{\bf A}]=
 {\cal D}_m\begin{pmatrix} Z^{n,k}  \\ Z^{n,k-1}\\ ...\\
 Z^{n,k-m+1}\\  Z^{n,k-m}
 \end{pmatrix}= \sum^{m}_{s=0}C^s_me^{{\i}(m-2s)\varphi}[{\cal D}Z^{n,k-s}]
 \nonumber \\
 &=& \frac {2e^{{\i}(m+n-2k)\varphi}}{n+1}\sum^{k}_{s=k-n}C^s_m
 \times
 \begin{cases}
 \cos[(n+1)(\beta-\varphi)] \  \ \text{for} \ n=\text{even} \\
 \\
 {\i}\sin[(n+1)(\beta-\varphi)] \  \text{for} \ n=\text{odd}\\
 \end{cases}\nonumber \\
 &=& \frac {2e^{{\i}(m+n-2k)\varphi}}{n+1}\alpha^{(m)}_{n,k}
 \times
 \begin{cases}
 \cos[(n+1)(\beta-\varphi)] \  \ \text{for} \ n=\text{even} \\
 \\
 {\i}\sin[(n+1)(\beta-\varphi)] \  \text{for} \ n=\text{odd,}\\
 \end{cases}\nonumber
\end{eqnarray}
 where $\alpha^{(m)}_{n,k}$
 are defined by the formula (\ref {alph1}).

 \noindent  Analogically, using the formula (\ref {Vtr2}) and Theorem 2
  we get
 \begin{eqnarray}
 && [{\cal D}_m{\bf B}]= {\cal D}_m\begin{pmatrix} \overline
 Z^{n,k-m} \\ \overline Z^{n,k-m+1}\\ ...\\ \overline Z^{n,k-1}\\ \overline
 Z^{n,k} \end{pmatrix}=\sum^{m}_{s=0}C^s_me^{{\i}(2s-m)\varphi}
 [{\cal D}\overline Z^{n,k-s}] \nonumber \\
 &=& \frac {2e^{{-\i}(m+n-2k)\varphi}}{n+1}\sum^{k}_{s=k-n}C^s_m
 \times
 \begin{cases}
 \cos[(n+1)(\beta-\varphi)] \  \ \ \ \text{for} \ n=\text{even} \\
 \\
 {-\i}\sin[(n+1)(\beta-\varphi)] \  \text{for} \ n=\text{odd}\\
 \end{cases}\nonumber \\
 &=& \frac {2e^{{-\i}(m+n-2k)\varphi}}{n+1}\alpha^{(m)}_{n,k}
 \times
 \begin{cases}
 \cos[(n+1)(\beta-\varphi)] \  \ \ \ \text{for} \ n=\text{even} \\
 \\
 {-\i}\sin[(n+1)(\beta-\varphi)] \  \text{for} \ n=\text{odd,}\\
 \end{cases}\nonumber
\end{eqnarray}

 From two formulas, derived above, it follows that
 \begin{eqnarray}
 &&\left [{\cal D}_m{\bf s}^{(+m)}_{n,k}
 \right ](\beta,\varphi)=
 \left [{\cal D}_m{\bf S}^{(+m)}_{n,k}
 \right ](\beta,\varphi)= (-1)^n[{\cal D}_m({\bf A}+{\bf B})] \nonumber  \\
 &=& \frac {2(-1)^n}{n+1}
 \alpha^{(m)}_{n,k}
 (e^{{\i}(n-2k+m)\varphi} \pm  e^{-{\i}(n-2k+m)\varphi}) \times
 \begin{cases}
  \cos[(n+1)(\beta-\varphi)] \\
  \\
 {\i}\sin[(n+1)(\beta-\varphi)]
 \end{cases} \nonumber  \\
 &=&\frac {4\alpha ^{(m)}_{n,k}}{n+1}
 \times
 \begin{cases}
 \cos [(n-2k+m)\varphi]\cos[(n+1)(\beta-\varphi)]  \\
 \\
 \sin  [(n-2k+m)\varphi] \sin [(n+1)(\beta-\varphi)]
 \end{cases}  \nonumber  \\
 &=&\frac {4||{\bf s}^{(+m)}_{n,k}||^2}{\pi 2^{m+1}}
 \times
  \begin{cases}
  \cos [(n-2k+m)\varphi] \cos[(n+1)(\beta-\varphi)] &\text {for \ n=even} \\
  \\
  \sin [(n-2k+m)\varphi] \sin[(n+1)(\beta-\varphi)] &\text {for \ n=odd}.
 \end{cases} \nonumber
 \end{eqnarray}
 By the same way we evaluate the fan-beam transform of the
 other part of basis for $n+m\neq 2k.$
 \begin{eqnarray}
 &&\left [{\cal D}_m{\bf s}^{(-m)}_{n,k}
 \right ](\beta,\varphi)=\left [{\cal D}_m{\bf S}^{(-m)}_{n,k}
 \right ](\beta,\varphi)=\dfrac{1}{{\i}}[{\cal D}_m({\bf A}-{\bf B})]
  \nonumber \\
 &&=\frac {2 \alpha^{(m)}_{n,k}}{n+1}
 \big (e^{{\i}(n-2k+m)\varphi} \mp e^{-{\i}(n-2k+m)\varphi} \big ) \times
 \begin{cases}
 \dfrac{1}{{\i}}\cos[(n+1)(\beta-\varphi)] \\
 \\
 \sin[(n+1)(\beta-\varphi)] \\
 \end{cases} \nonumber  \\
 &&=\frac {4\alpha^{(m)}_{n,k}}{n+1}
 \times
  \begin{cases}
  \sin [(n-2k+m)\varphi] \cos[(n+1)(\beta-\varphi)] \\
  \\
  \cos [(n-2k+m)\varphi] \sin[(n+1)(\beta-\varphi)] \\
 \end{cases}   \nonumber \\
 &&=\frac {4||
 {\bf s}^{(-m)}_{n,k}||^2}{\pi 2^{m+1}}
 \times
\begin{cases}
  \sin [(n-2k+m)\varphi] \cos[(n+1)(\beta-\varphi)]  &\text {for \ n=even} \\
  \\
  \cos [(n-2k+m)\varphi] \sin[(n+1)(\beta-\varphi)]  &\text {for \ n=odd}. \\
\end{cases} \nonumber
  \end{eqnarray}
 \\
 Consider the case  $n+m=2k$ and $n, \ m=$ even, then
 \begin{eqnarray}
 \left [{\cal D}_m{\bf s}^{(+m)}_{n,\frac {m+n}{2}}
 \right ](\beta,\varphi)=
 \left [{\cal D}_m{\bf S}^{(+m)}_{n,\frac {m+n}{2}}
 \right ](\beta,\varphi)=
 \frac {4\alpha ^{(m)}_{n,k}}{2(n+1)}
 \times
 \begin{cases}
 \cos[(n+1)(\beta-\varphi)]   \\
 0
 \end{cases} \nonumber
 \end{eqnarray}
 \begin{eqnarray}
 =\frac {4||{\bf s}^{(+m)}_{n,\frac {m+n}{2}}||^2}{\pi 2^{m+1}}
 \times
 \begin{cases}
 \cos[(n+1)(\beta-\varphi)]  &\text {for \ n=even} \\
 0                           &\text {for \ n=odd}.
 \end{cases} \nonumber
 \end{eqnarray}
 If  $n+m=2k$ and  $n, \ m=$  odd, then
 \begin{eqnarray}
 \left [{\cal D}_m{\bf s}^{(-m)}_{n,\frac {m+n}{2}}
 \right ](\beta,\varphi)=
 \left [{\cal D}_m{\bf S}^{(-m)}_{n,\frac {m+n}{2}}
 \right ](\beta,\varphi)= \frac
 {4\alpha^{(m)}_{n,k}}{2(n+1)}\times
 \begin{cases}
  0  \\
  \sin[(n+1)(\beta-\varphi)]
 \end{cases}  \nonumber
 \end{eqnarray}
\begin{eqnarray}
 =\frac {4||{\bf s}^{(-m)}_{n,\frac {m+n}{2}}||^2}{\pi
 2^{m+1}}\times
 \begin{cases}
  0                           &\text {for \ n=even}  \\
 \sin[(n+1)(\beta-\varphi)]   &\text {for \ n=odd.}
 \end{cases} & \nonumber
 \end{eqnarray}
 Using equations for norms (\ref{alph}), we get  (\ref{Da}).
 Theorem 3 is proved.
 \endproof
 \par At the end of this section  we present some examples.
 \par \example 1 Let's take for instance $m=0,$ that corresponds to the scalar field,
 hence  we have  $a(x^1,x^2)=A(z,\overline z)$ and  the orthogonal basis
  ${s}^{(\pm 0)}_{n,k}={S}^{(\pm 0)}_{n,k}$ in  $L_2({\mathbb D})$ is
 \begin{eqnarray}
 \label{m0+}
 {s}^{(+0)}_{n,k}&=&
 \begin{cases} (-1)^n 2 \Re \ Z^{n,k} \  &\text {for \ $2k \neq n$}   \\
 & \nonumber \\
 s^{(+0)}_{n,\frac {n}{2}}=  Z^{n,\frac {n}{2}} \ &\text {for \ $2k = n,$}
 \end{cases}
 \nonumber
 \\
 && \nonumber \\
  \label{m0-}
  {s}^{(-0)}_{n,k}&=&
 \begin{cases}
 2\Im \ Z^{n,k} \ \ \ \ \ \ \ \              &\text { for \ $2k \neq n$} \\
 &\\
 s^{(-0)}_{n,\frac {n}{2}}=0 \ \ \ \ \ \ \ \ &\text { for \ $2k = n,$}  \\
 \end{cases}
 \nonumber
 \end{eqnarray}
 where $n \geq 0, \ k=0,1,...,\left [ \frac {n}{2} \right ].$
 We have
 $\dim {\bf H}_N(\mathbb D;{\bf S}_0, \delta=0)= \frac{(N+1)(N+2)}{2}$
 and singular values  are
 $$\sigma_{n,k}\equiv\sigma^{(\pm 0)}_{n,k}=
  \sqrt{\frac {8\pi}{n+1}}, \
  \left ( n \geq 0,  \  k=0,...,\left [\frac {n}{2}\right ] \right ).
 $$
 \par \example 2
  For $m=1$ one gets covector field
  ${\bf a}(x^1,x^2)=\{a_{1}, a_{2}\},$
  which in the complex variables according to the tensor law has
  the representation
  ${\bf A} (z,\overline z)=
  \left \{ \dfrac {a_1-{\i}a_2}{2}, \dfrac {a_1+{\i}a_2}{2} \right \}.$
 Dimension of the finite-dimensional subspace ${\bf H}_N({\mathbb D;\bf S}_1,\delta=0)$
 equals to $\frac{(N+1)(N+4)}{2}$ and
 singular values  are
  $$\sigma^{(\pm 1)}_{n,k}=
  \begin{cases}
  {\sqrt {\dfrac {4\pi}{n+1}} \  \ \text{if} \ n \geq 0 \ \text{and}   \ k=0} \\
  {\sqrt {\dfrac {8\pi}{n+1}} \ \ \text{if} \ n \geq 1 \ \text{and}
  \ k=1,...,\left [\dfrac {n+1}{2}\right ]}.
  \end{cases}
  $$
  The polynomial orthogonal basis ${\bf s}^{(\pm 1)}_{n,k}$ of the space
  of solenoidal covectors fields ${\bf H}(\mathbb D;{\bf S}_1, \delta=0)$ looks as
  follows
 \begin{eqnarray}
 \label{m1+}
 {\bf s}^{(+1)}_{n,k} \rightarrowtail
 {\bf S}^{(+1)}_{n,k}&=&
 \begin{cases}(-1)^n
 \left \{ Z^{n,k} + {\overline {Z}}^{n,k-1},Z^{n,k-1} + {\overline {Z}}^{n,k}  \right \}
 &\text{for} \ 2k \neq n+1 \nonumber   \\
 & \ \nonumber \\
 {\bf S}^{(+1)}_{n,\frac {n+1}{2}}= \{0,0 \}
 \  &\text {for} \ 2k =n+1,
 \nonumber
 \\
 \end{cases}
 \nonumber \\
 && \nonumber\\
 \label{m1-}
 {\bf s}^{(-1)}_{n,k} \rightarrowtail
 {\bf S}^{(-1)}_{n,k}&=&
 \begin{cases}
 \dfrac{1}{{\i}}  \left \{Z^{n,k} - {\overline {Z}}^{n,k-1},Z^{n,k-1} -{\overline {Z}}^{n,k}  \right \}
  \ \ \ \ \ \ &\text {for} \ 2k \neq n+1  \nonumber  \\
 & \nonumber \\
 {\bf S}^{(-1)}_{n,\frac {n+1}{2}}=\dfrac{1}{{\i}}
 \left  \{ Z^{n,\frac {n+1}{2}}, Z^{n,\frac {n+1}{2}-1}  \right \}
  \ \ \ \ \ \ &\text {for} \  2k =n+1, \\
 \end{cases}
 \nonumber
 \end{eqnarray}
 where $\ n \geq 0, \ k=0,1,...,\left [ \frac {n+1}{2} \right ].$
 \par \example 3
  For $m=2$  we have a symmetric second-order 2D tensor field
  $ {\bf a}(x^1,x^2)=
  \left \{\begin{array}{ll}
  a_{11} & a_{12}\\
  a_{21} & a_{22}\\
  \end{array} \right \},
  $
 which in complex variables has components
 $${\bf A}(z,\overline z)= \left \{
 \begin{array}{cc}
 A_{11} & A_{12}\\
 A_{21} & A_{22}\\
 \end{array}
 \right \}
  = \left \{
 \begin{array}{ll}
 \dfrac {a_{11}-a_{22}-2{\i}a_{12}}{4} &\dfrac {a_{11}+a_{22}}{4}\\
 \dfrac {a_{11}+a_{22}}{4} & \dfrac {a_{11}-a_{22}+2{\i}a_{12}}{4}\\
 \end{array}
 \right \}.
 $$
  We also  have inverse equalities
 $$ {\bf a}(x^1,x^2)=
  \left \{\begin{array}{ll}
 2(A_{12}+\Re A_{11}) & -2 \Im A_{11}\\
 -2 \Im A_{11}       & 2(A_{12}-\Re A_{11})\\
  \end{array} \right \}.
 $$
 Singular values for this case are
 $$\sigma^{(\pm 2)}_{n,k}=
 \begin{cases}
 {\sqrt  {\dfrac {2\pi}{n+1}} \ \ \text{if} \ n \geq 0 \ \text{and}   \ k=0} \\
 {\sqrt {\dfrac {4\pi}{n+1}} \  \ \text{if} \ n =    0 \ \text{and}   \ k=1} \\
 {\sqrt {\dfrac {6\pi}{n+1}} \ \ \text{if} \ n \geq 1 \ \text{and}   \ k=1} \\
 {\sqrt {\dfrac {8\pi}{n+1}} \ \ \text{if} \ n \geq 2 \ \text{and}   \
  k=2,...,\left [\frac {n+2}{2}\right ]},
 \end{cases}
 $$
 where $n \geq 0, k=0,1,...,\left [ \frac {n+2}{2} \right ]$
 and basis tensor fields are
 \begin{eqnarray}
 \label{m2+}{\bf s}^{(+2)}_{n,k} \rightarrowtail
{\bf S}^{(+2)}_{n,k}&=&
  \begin{cases}
  & (-1)^n \left \{\begin{array}{ll}
  Z^{n,k} + \overline {Z}^{n,k-2}   & Z^{n,k-1}+\overline {Z}^{n,k-1}  \\
  Z^{n,k-1} + \overline {Z}^{n,k-1} & Z^{n,k-2}+ \overline {Z}^{n,k} \\
  \end{array} \right \}  \\
  &\text{for}  \ 2k \neq n+2  \\
  & \\
  &{\bf S}^{(+2)}_{n,\frac {n+2}{2}}=
  \left \{\begin{array}{ll}
  Z^{n,\frac {n+2}{2}} & Z^{n,\frac {n+2}{2}-1}  \\
  Z^{n,\frac {n+2}{2}-1}  & Z^{n,\frac {n+2}{2}-2} \\
  \end{array} \right \}  \\
  &\text{for}  \ 2k =n+2,
  \end{cases}   \nonumber \\
 && \nonumber\\
 \label{m2-}
 {\bf s}^{(-2)}_{n,k} \rightarrowtail
 {\bf S}^{(-2)}_{n,k}&=&
  \begin{cases}
  &\dfrac {1}{{\i}}
  \left \{ \begin{array}{ll}
  Z^{n,k} - \overline {Z}^{n,k-2} & Z^{n,k-1}-\overline {Z}^{n,k-1}  \\
  Z^{n,k-1} -\overline {Z}^{n,k-1} & Z^{n,k-2}- \overline {Z}^{n,k} \\
  \end{array} \right \} \\
  &\text {for} \ 2k \neq n+2 \\
  &\\
  &{\bf S}^{(-2)}_{n,\frac {n+2}{2}}= \dfrac {1}{{\i}}
  \left \{ \begin{array}{ll}
  0 & 0  \\
  0 & 0 \\
  \end{array} \right \} \\
  &\text {for} \ 2k =n+2.
 \end{cases}
 \nonumber
\end{eqnarray}
Also in this case we have $\dim {\bf H}_N(\mathbb D;{\bf S}_2, \delta=0)= \frac{(N+1)(N+6)}{2}.$

\section{Implementation}
 Scalar and vector cases of the inversion formula were numerically
 implemented and tested. The algorithm consists of 3 parts:
 solving the direct problem (that emulates the data acquisition in real life),
 finding coefficients of the polynomial that represents
 a function being reconstructed and evaluation of this polynomial on a
 grid for visualization.
 \par
 In the scalar case, given a test function, defined by its values on a rectangular grid,
 the direct problem was solved by
 computing integrals (\ref {scfanT}) for the number of discrete values
 $\beta_p,$ $\varphi_q$
 \begin{equation}
 [Da](\beta_p,\beta_p-\frac{\pi}{2}+\varphi_q)=f_{p,q}, \ (p,q=0,1,...,M+1).
 \label{i1}
 \end{equation}
 \\
 Bilinear interpolation was used to get
 the values of the original function between knots. So the obtained
 data set is an $(M+2)\times(M+2)$ matrix of $(f_{p,q})$ values
 that serves as an input for the
 inversion algorithm. Consider the scalar case for instance. Then,
 the function $a(x,y)$ is approximated by the
 polynomial of degree $N\leq M$
 (note, that in this section we use notations $x\equiv  x^1$
 and $y\equiv x^2$)
 \begin{eqnarray}
 a_N(x,y)&=&2\sum_{n=0}^{N}\sum^{[n/2]\ast}_{k=0}
 \Bigl[
 \bigl(a_{n,k}\cos((n-2k)\psi)-b_{n,k}\sin((n-2k)\psi)\bigr) \nonumber \\
 &&\times(-1)^k r^{n-2k}P_k^{(0,n-2k)}(2r^2-1) \Bigr].
 \label {i2}
 \end{eqnarray}
 \par
 Here $(r,\psi)$ are the polar coordinates of the point $(x,y)$
 and the sign $\ast$ means that in the case of even $n$ the coefficient
 $a_{n,[n/2]}$ should be divided by $2$ and $b_{n,[n/2]}$ should
 be set to $0.$
 Then the fan-beam transform of (\ref{i2}) will look like
 \begin{equation}
 [{\cal D}a_N](\beta,\varphi)=\sum_{n=0}^{N}\frac{4}{n+1}
 \sum_{k=0}^{[n/2]\ast}a_{n,k}\times
 \begin{cases}
 \cos[(n+1)(\beta-\varphi)]\cos(n-2k)\varphi \\
 \sin[(n+1)(\varphi-\beta)]\sin(n-2k)\varphi
 \end{cases}$$
 $$-b_{n,k}\times
 \begin{cases}
 \cos[(n+1)(\beta-\varphi)]\sin(n-2k)\varphi \\
 \sin[(n+1)(\beta-\varphi)]\cos(n-2k)\varphi,
 \end{cases}
 \label{i3}
 \end{equation}
 where the upper lines in the braces are used for the even $n$ and the
 lower lines --- for the odd $n.$ The sign $\ast$ means the same
 as in (\ref{i2}). After the substitution of (\ref{i3}) into the
 (\ref{i1}) we get a system of linear
 equations for determining  $\dfrac {(N+1)(N+2)}{2}$ unknown
 coefficients $a_{n,k}$ and $b_{n,k}.$
 \par
 In the case of regular scanning scheme $\beta_p=p\varepsilon,$
 $\varphi_q=q\dfrac{\varepsilon}{2},$
 $\varepsilon=\dfrac{2\pi}{M+2}$ an explicit formulas for
 determining coefficients $a_{n,k}$ and $b_{n,k}$ were derived,
 provided that $M=N$
 \begin{eqnarray}
 a_{n,k}=(-1)^{k}\frac{n+1}{(M+2)^2}
 \sum_{p=0}^{M+1}\sum_{q=0}^{M+1}f_{p,q}
 \sin\left[\varepsilon\left(p(2k-n)+\frac{q}{2}(2k+1)\right)\right],
 \label{i6}  \end{eqnarray}
 \begin{eqnarray}
 b_{n,k}=(-1)^{k+1}\frac{n+1}{(M+2)^2}
 \sum_{p=0}^{M+1}\sum_{q=0}^{M+1}f_{p,q}
 \cos\left[\varepsilon\left(p(2k-n)+\frac{q}{2}(2k+1)\right)\right].
  \label{i7}
 \end{eqnarray}
 Analogical formulae for the parallel-beam geometry can be found in \cite{Mar74}.
 The implementation of formulas (\ref{i6}) and (\ref{i7}) uses  FFT and requires
 $\mathcal{O}(N^2\log_2N)$ operations, see \cite{BK01}.

 After the coefficients $a_{n,k}$ and $b_{n,k}$ are found,
 the polynomial that represents the reconstructed function is
 effectively evaluated using a recurrent formula, see \cite{BK01}.

 The inversion algorithm was also tested under the presence of
 a noise in the input data (sinogram). Uniform and Poisson random
 distributions were used for this purpose.

 A representative set of numerical tests was performed for scalar
 and vector cases of the inversion algorithm. Some of the results are
 shown on figures 2--4.

 On the top of the figure 2 there are original function (to be
 reconstructed) on the left-hand side and its sinogram
 (an input data set for the inversion algorithm) on the right-hand side.
 The middle row contains reconstructions from 32 and 256 fan-beam
 projections (free of noise). The number of terms in SVD was 30
 and 254 respectively. The bottom row contains examples of
 reconstruction from noisy data. A random noise was superimposed on
 the sinogram. The $L_2$-norm of the noise was 10\% of the
 $L_2$-norm of the sinogram. 1024 noisy fan-beam projections were
 used. The number of terms in SVD that were taken for
 reconstruction are 1022 and 254 respectively. It's possible to
 reduce the noise in the output image by taking less terms in the
 SVD.

 Another example of scalar tomography is shown on the figure 3.
 Again, the original unknown function (the fast oscillating one,
 with fine features) is at the top row, on the
 left-hand side and it's sinogram is on the right-hand side.
 The middle row contains reconstructions from 32 and 512 fan-beam
 projections (free of noise). The number of terms in SVD was 30
 and 510 respectively. The bottom row contains examples of
 reconstruction from noisy data. A random noise was superimposed on
 the sinogram. The$L_2$-norm of the noise was 10\% of the
$L_2$-norm of the sinogram. 2048 noisy fan-beam projections were
 used. The number of terms in SVD that were taken for
 reconstruction are 1022 and 510 respectively. Again, one can observe significant
 enhancement of reconstruction when only part of terms are taken in SVD.

 The figure 4 illustrates the vector case of the inversion algorithm. The
 first solenoidal vector field (the top row, where the first component
  $a_1$ is on the left
 and the second component $a_2$ is on the right) is defined by the formulae
 \begin{eqnarray}
 a_1(x,y)&=&2xy\cos(x^2+y^2)+\cos(6xy)-6xy\sin(6xy),\nonumber\\
 a_2(x,y)&=&-\sin(x^2+y^2)-2x^2\cos(x^2+y^2)+6y^2\sin(6xy).
 \label{i8}
 \end{eqnarray}
  Another vector field (the middle row) was obtained from the previous
  solenoidal vector field by adding the potential part
 \begin{align}
 a_1(x,y) & =2xy\cos(x^2+y^2)+\cos(6xy) -6xy\sin(6xy)
 \nonumber \\
 &+2\pi x\cos(\pi(x^2+y^2)), \nonumber \\
 a_2(x,y) & =-\sin(x^2+y^2)-2x^2\cos(x^2+y^2) +6y^2\sin(6xy)\nonumber \\
 &+2\pi y\cos(\pi(x^2+y^2)).
 \label{i9}
 \end{align}
 As it can be seen (on the bottom row of the figure), reconstruction
 from these two vector fields is identical and contains only the
 solenoidal part of the vector field.

 The figure 5 illustrates another solenoidal vector field (the top
 row). Its reconstruction from the 20 irregular fan-projections is in the middle
 row.   Here the positions of the fan-projection centers are shown
 as white dots on the boundary of the circle.
 Scanning was performed only over those lines, whose endpoints
 belong to this set of 20 points.
 The last reconstruction (the bottom row) was made under the presence of
 noise in the sinogram. The noise level was 3\% (again, in $L_2$-norm).

\begin{figure}[p!]
 \begin{center}
 \includegraphics[scale=0.55, bb= 0 0 453 680]{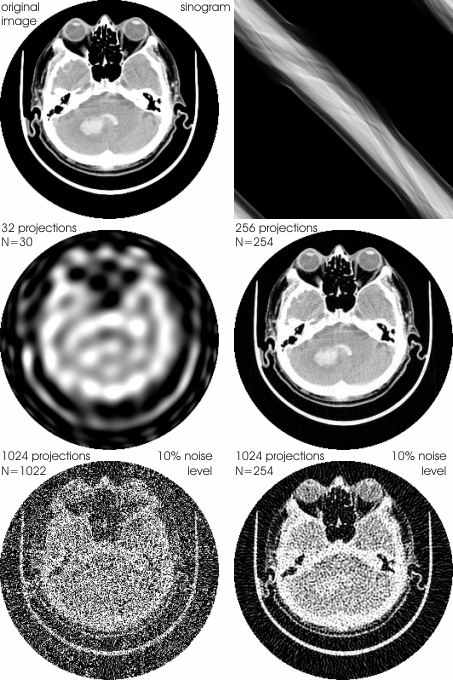}
 \caption{\footnotesize Top row: image of the original function (on the left)
  and its sinogram (on the right).
 Middle and bottom rows: reconstructions of the function from different number
 of fan-projections and under the
 presence of noise in the sinogram, see the text for detailed explanation.}
 \end{center}
\end{figure}

 \begin{figure}[p!]
 \begin{center}
 \includegraphics [scale=0.55, bb=0 0 453 680]{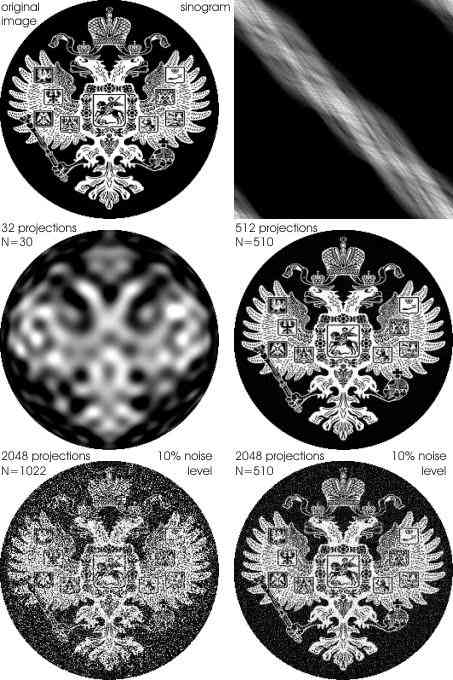}
 \caption{\footnotesize Top row: image of the original function (on the left) and its
  sinogram (on the right).
 Middle and bottom rows: reconstructions of the function from different number
 of fan-projections and under the
 presence of noise in the sinogram, see the text for detailed explanation.}
 \end{center}
\end{figure}

\begin{figure}[p!]
 \begin{center}
 \includegraphics [scale=0.55, bb=0 0 453 679]{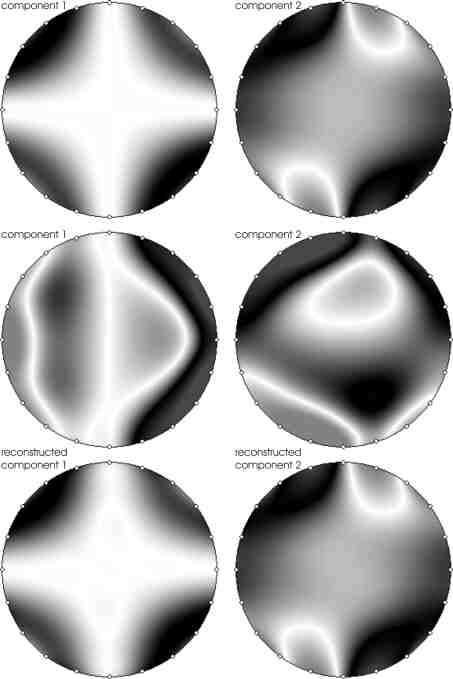}
 \caption{\footnotesize Top row: solenoidal vector field (first component on the left,
  second component on the right)
 given by formulae (\ref{i8}). Middle row: non-solenoidal vector field given by formulae
  (\ref{i9}). Bottom row:
 reconstruction from 20 fan-projections is identical in both cases, relative
 error $0.21\%.$}
 \end{center}
\end{figure}

 \begin{figure}[p!]
 \begin{center}
 \includegraphics [scale=0.55, bb =0 0 453 679] {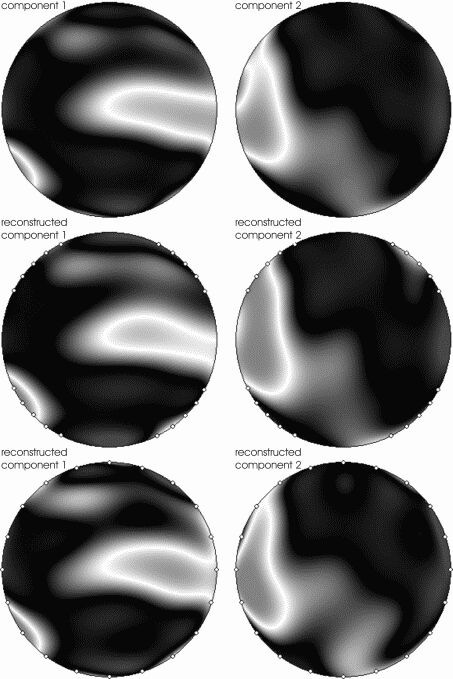}
 \caption{\footnotesize Top row: original solenoidal vector field. Middle row:
  its reconstruction from 20 irregular
 fan-projections, relative error $2.6\%.$
 Bottom row: its reconstruction from the 20 regular
 fan-projections under the presence of noise ($3\%$)
 in the sinogram, relative error $1.6\%.$}
 \end{center}
\end{figure}

\section{Conclusion}

 The novel inversion algorithm for the tensor tomography problem was
 developed and numerically implemented. The algorithm is based on the SVD of the
 tensor Radon transform that allows to characterize the range of the operator,
 invert it and estimate an incorrectness of the corresponding inverse problem.
 The algorithm can also be used with noisy measurements.

\refs

 \bibitem{ABK98}
 E.V.~Arbuzov, A.L.~Bukhgeim, and S.G.~Kazantsev,
 Two-dimentional tomography problems and the theory of A-analytic function.
 {\refstyle Siberian Advances in Mathematics\/} (1998) \vol{8}\no 4, 1--20.

 \bibitem{Bo01} J.~Boman, \  Injectivity for a weighted vectorial Radon transform,
 {\refstyle Contemp. Math.\/} (2001), \vol{278}, 87--95.

 \bibitem{BK01}
 A.A.~Bukhgeim and S.G.~Kazantsev,
 {\it Singular value decomposition of the 2D fan-beam Radon transform of
 tensor fields in a unit disc.}
 Novosibirsk, 2001. (Preprint/RAS. Siberian Dep.  Institute of
 Mathematics\no 86) (in Russian).


 \bibitem {BW83}
  M.~Born and  E.~Wolf,
 {\it  Principles of Optics.}
  Pergamon, New-York,  1983.

 \bibitem{BH}
  H.~Braun  and A.~Hauck,
 Tomographic reconstruction of vector fields.
 {\refstyle IEEE Transaction on Signal  Processing\/} (1991)
 \vol{39}\no 2, 464--471.

 \bibitem {Corm63}
 A.M.~Cormack, Representation of a function by its
 line integrals, with some radiological applications. I
 {\refstyle J. Appl. Phys.\/} (1963)\no 34, 2722--2727.

 \bibitem{Corm64}
 A.M.~Cormack, Representation of a function by its line
 integrals with some radiological applications.  II
 {\refstyle J. Appl. Phys.\/} (1964)\no 35, 195--207.

 \bibitem{DL}
 R.~Dautray and J.-~L.~Lions,
 {\it Mathematical Analysis and Numerical Methods for Science and Technology.}
 V.3. Spectral Theory and  Applications. Springer, 2000.

 \bibitem{BDS02}
 M.A.~Bezuglova, E.Yu.~Derevtsov, and S.B.~Sorokin,
 The reconstruction of a vector field by finite difference methods.
 {\refstyle J. Inverse Ill-posed Problems\/} (2002) \vol{10}\no 2, 125--154.

 \bibitem{DK02}
 E.Yu.~Derevtsov and I.G.~Kashina,
 Numerical solution to the vector
 tomography problem by tools of polynomial basis.
 {\refstyle Sib. J. Numerical Math.} (2002) \vol{5}\no 3, 233--254 (in Russian).

 \bibitem{Ca81}
 M.~Cantor,  Elliptic operators and the decomposition of tensor fields.
 {\refstyle Bull. AMS\/} (1981) \vol{\bf 5}\no 3, 1981,
 235-262.


 \bibitem{Herm80}
  G.T.~Herman,
  {\it Image Reconstruction from
 Projections: The Fundamentals of Computerized Tomography.}
 Academic Press, New-York, 1980.

 \bibitem{GR80}
 I.S.~Gradshteyn and  I.M.~Ryzhik,
 {\it Tables of integrals, sums, series and products.}
 Academic Press,  New-York, 1965.

 \bibitem{G-MA94}
 O.~Gil-Medrano and  A.M.~Amilibia,
 About a Decomposition of the Space of Symmetric Tensors
 of Compact Support on a Riemannian Manifold.
 {\refstyle New York J. Math.\/} {1994} \vol{1}, 10--25.

 \bibitem{M90}
 P.~Maass,
 Singular value decompositions for Radon transform.
 In: {\it Mathematical Mathods in Tomography}.
 Springer-Verlag, 1990, 6--14.

 \bibitem{Mar74}
 R.B.~Marr,
 On the Reconstruction of a Function on a
 Circular Domain from a Sampling of its Line Integrals.
 {\refstyle J. of Mathematical Analysis and Applications\/}
 (1974) \vol{45}\no 2, 357--374.

 \bibitem {Nat86} F.~Natterer,
 {\it  The mathematics of  Computerized Tomography.}
  Teubner Verlag, Stuttgart, 1986.

 \bibitem{NW01}
 F.~Natterer  and F.~Wubbeling,
 {\it Mathematical Methods in Image Reconstruction.}
 Philadelphia, SIAM, 2001.

 \bibitem{Nor}
 S.J.~Norton,
 Tomographic reconstruction of 2-D vector fields:
 application to flow imaging.
 {\refstyle J. of Geophysics\/} (1987) \vol{97}, 161--168.


 \bibitem{PR89}
 A.~Prata and W.V.~Rusch,
 Algorithm for computation of Zernike polynomials expansion
 coefficient. {\refstyle Applied Optics\/} (1989) \vol{28}\no 4, 749--754.

\bibitem{Qui83}
 E.T.~Quinto,
 Singular Value Decomposition and Inversion Methods for the Exterior
 Radon Transform and thr Spherical Transform.
 {\refstyle J. of mathematical analysis and applications\/} (1983) \vol{95}\no 2, 437--448.

\bibitem{Sch01}
 Th.~Schuster,
 An efficient mollifier method for three-dimensional
 vector tomography: convergence analysis and implementation.
 {\refstyle Inverse Problems\/} (2001) \vol{17}, 739--766.

 \bibitem{Shar}
 V.A.~Sharafutdinov, {\it Integral Geometry of Tensor Fields.}
 VSP, Utrecht, 1994.


 \bibitem{St99}
 K.~Strahlen,
 Studies of Vector Tomography.
 {\it Ph.D. thesis. } Lund University, 1999.

 \bibitem{Ver00}
 L.V.~Vertgeim,
 Integral geometry problems for symmetric tensor fields with incomplete data.
 {\refstyle J.Inv. Ill-Posed Problems\/} (2000) \vol{8}\no 3,  353--362.
 \bibitem{V}
 I.I.~Vekua, {\it Tensor analysis.} Nauka, 1988.

 \bibitem{W}
 H.~Weyl, The method of orthogonal projection in potential theory.
 {\refstyle Duke Math. J.\/} (1940)\no 7, 411--444.
\endrefs

\end{document}